\documentclass[12pt]{article}
\usepackage{amssymb,latexsym,amsmath}
\def\C{\mathbb C}
\def\R{\mathbb R}
\def\N{\mathbb N}
\def\Z{\mathbb Z}

\newtheorem{thm}{Theorem}[section]
\newtheorem{lem}{Lemma}[section]
\newtheorem{con}{Conjecture}[section]
\newtheorem{definition}{Definition}[section]
\newtheorem{prop}{Proposition}[section]

\begin{document}
\sffamily
%\maketitle
%\thispagestyle{empty}

\title{A special case of a conjecture of \\ Hellerstein, Shen and Williamson}
% and differential polynomials}
\author{J.K. Langley}
\maketitle

\begin{center}
\emph{Dedicated to the memory of I.V. Ostrovskii}
\end{center}

\begin{abstract}
The paper proves a special case of a conjecture of Hellerstein, Shen and Williamson concerning non-real zeros of 
derivatives of real meromorphic
functions. \\
Keywords: meromorphic function, non-real zeros. 
MSC 2010: 30D20, 30D35.
\end{abstract}

%CHECKEDANDREVISEDMARCH23NOTUPDATEDONARXIV

\section{Introduction}

This paper concerns
%orem \ref{thmrealmero} represents a very special case of
 the problem of classifying real
meromorphic functions in the plane which, together with some of their derivatives,
have only real zeros and poles 
\cite{HSW,Hin1,Hin3}.  Here a meromorphic function $f$ is called real if   $f(\R) \subseteq \R \cup \{ \infty \}$, and it is known that if $f$ is 
real entire and $f$ and $f''$ have only real zeros  then $f$ 
belongs to the Laguerre-Polya class $LP$, as conjectured by  Wiman and proved in 
\cite{BEL,LeO,SS}.  For the meromorphic case,  the following conjecture was advanced in
\cite{HSW}.

\begin{con}[\cite{HSW}]
 \label{conhsw}
Let $f$ be a real transcendental meromorphic function in the plane with at least one pole, and assume that all zeros and poles of
$f$, $f'$ and $f''$ are real, and that all poles of $f$ are simple. Then 
$f$ satisfies
%there exist real numbers $A$, $B$, $C$, $a$, $b$ such that
\begin{equation}
 \label{hswform}
f(z) =  C  \tan (az+b) + Dz + E, \quad  a, b, C, D, E \in \R. 
\end{equation}
\end{con}
In the absence of the assumption that  $f$ has only simple poles, further examples arise
for which $f$, $f'$ and $f''$ have only real zeros and poles
 \cite{Hin4}.  Conjecture \ref{conhsw} is known to be true if any of the following additional hypotheses holds:\\
(a)  $f'$ omits some finite value \cite{HSW,HinRos,Laams09,ankara,Lawiman13,nicks,rossireal};\\
(b)  $f$ has infinitely many poles and $f''/f'$ has finitely many zeros \cite[Theorem 1.5]{Lajda18};\\
(c) $f$ has infinitely many zeros and poles, all real, simple and \textit{interlaced} --  that is,  
between any two consecutive poles of $f$ there is a zero, and 
between consecutive zeros of $f$ lies a pole \cite{wiman2022}.

%It was also shown in  \cite{Lawiman10,wiman2022} that there are no functions $f$
%which satisfy the hypotheses of Conjecture~\ref{conhsw} as well as one of the following two conditions:
%$f$ has finite order and infinitely many poles but finitely many zeros; 
%$f$ has infinite order and the zeros or poles of $f$ have finite exponent of convergence.

The following theorem implies a further special case of Conjecture \ref{conhsw}.  

\begin{thm}
\label{thm1}
Let $f$ be a real  meromorphic function in the plane, 
such that $f$ and $f'$ have no zeros or poles in $\C \setminus \R$, while  $f''/f$ has no zeros in $\C$.
Then
there exist $\alpha_1, \alpha_2, \alpha_3 \in \R$, with $\alpha_1 \alpha_2  \neq 0$,
such that 
$g(z) = \alpha_1 f(\alpha_2 z + \alpha _3 )$ satisfies one of the following:\\
(i) $g(z) = f_1(z) = \sin z$;\\
(ii) $g(z) = f_2(z) = e^z $;\\
(iii) $g(z) = f_3(z) = \tan z $;\\
(iv) $g(z) = f_4(z)$, where
\begin{equation}
 f_4(z) = \sum_{k=0}^\infty \frac{z^{k+1}}{k!(k+1)!} = 
z + \frac{z^2}2 + \frac{z^3}{12} + \ldots 
\label{besselsoln}
\end{equation}
solves
\begin{equation}
z y''(z) =  y(z);
\label{bessel}
\end{equation}
(v)  $g(z) = F_1(z) = z^Q$ for some $Q \in \Z \setminus \{ 0, 1 \}$;\\
(vi) $g(z) = F_2(z)$, where 
\begin{equation}
 F_2(z) = 
 \frac{d^{n-2} }{dz^{n-2}} \left( z^{n-1} (z-1)^{n-1} \right)
  \label{hyper1afn} 
  \end{equation}
 for some integer $n \geq 2$, and $F_2$ solves 
 \begin{equation}
\label{hypergeometricA}
z (z-1)  y''(z) = n(n-1)  y(z) ;
\end{equation} 
(vii) $g(z) = F_3(z)$, where $F_3$ is given by 
\begin{equation} 
F_3(z) = 
(z-K) H_n \left(  \frac{K + 1}{K-1} -\frac{2K}{(K-1)z } \right) 
\label{hyper2Afn} 
\end{equation}
for some integer  $n \geq 1$ and $K \in \R \setminus \{ 0, 1 \}$, in which  
\begin{equation} 
\label{Hndefn}
H_n(w) =  \frac{d^n}{dw^n} \left( (w-1)^{n-1}(w+1)^{n+1}  \right) ,
\end{equation}
while $F_3$ solves 
\begin{equation}
\label{hypergeometric3A}
z^2 (z-1)(z-K)  y''(z) = Kn(n+1) y(z);
\end{equation} 
(viii) $g(z) = F_4(z)$, where 
\begin{equation}
F_4(z) = H_n \left(  1 -\frac2{z } \right) ,
\end{equation}
in which   $1 \leq n \in \N$ and $H_n$ is given by (\ref{Hndefn}), while $F_4$ solves 
\begin{equation}
\label{hypergeometric3AB}
z^2 (z-1) y''(z) = - n(n+1) y(z).
\end{equation}  
\end{thm} 

 Conversely, the equations (\ref{bessel}), (\ref{hypergeometricA}), (\ref{hypergeometric3A})   (for $K > 1$) and 
 (\ref{hypergeometric3AB}) all supply examples satisfying 
  the hypotheses of the theorem. 
The function $f_4$ in  (\ref{besselsoln}) and its connection to Bessel functions will 
be discussed  in Section \ref{besselexample}, while the 
rational functions $F_2, F_3, F_4$ in (vi), (vii) and (viii), which are linked to hypergeometric functions,
 will be treated in detail in Sections \ref{hyper1a},
\ref{hyper2A} and 
\ref{hyper2AA}.

Of course, the condition that $f''/f$ has no zeros makes 
Theorem \ref{thm1}  a very special case of Conjecture \ref{conhsw}, albeit without the assumption that $f$ is transcendental and all poles of $f$ are simple,
but the fact that the proofs of all the resolved special cases are   lengthy tends to suggest that the full conjecture is difficult. 
The result may also be viewed as a special case of the problem of 
determining 
all meromorphic functions $f$
such that $f''/f$ has no zeros: in this direction,   it was proved in \cite{jklcv} that if $f$ is entire of order less than $1$, or meromorphic of order less than $1/2$, and $f''/f$ is transcendental, then $f''/f$ has infinitely many zeros. 

Note that the corresponding problem for the case where $f$ is strictly non-real, that is,  $f$ is not a constant multiple of a real meromorphic function, was completely settled in \cite{HSWsnr}, the main result of which classified all strictly non-real meromorphic functions $f$ in the plane for which $f$, $f'$ and $f''$ have only real zeros and poles.  
%, as well as in  Proposition \ref{proprational}. 

In common with much of the work on non-real zeros of derivatives, this paper relies heavily on key results and methods developed by B. Ja. Levin and I.V. Ostrovskii, as set out in the paper \cite{LeO} and the textbook \cite{GO} - in particular, the factorisation of the logarithmic derivative 
(Section~\ref{levostfact}) and an integral inequality  linking  the Tsuji and Nevanlinna proximity functions 
(Lemma~\ref{lemlevinost}).

\section{Preliminaries and examples}\label{examples}

First, let $D$ be a real entire function, whose zeros $x_k$ are all real and simple. Then a standard application
of the Mittag-Leffler theorem gives a real entire function $C$ such that 
$$
\frac{e^{C(z)}}{D(z)} = \frac{-2}{z-x_k} + O( |z-x_k|) 
$$
as $z \to x_k$, for each $k$. The formula  $g'/g = e^C/D$  then defines
 a real meromorphic function
$g$, 
such that $g$ and $g'$ have no zeros at all, while  for each $k$ there exists $c_k \in \R \setminus \{ 0 \}$
with $g(z) = c_k (z - x_k)^{-2} + O(1)$ as $z \to x_k$. Hence there exists a real   meromorphic function $f$
with $f' = g$ and $f'f''$ zero-free. However, this construction of course gives no control over the
location of the zeros of $f$ itself. 

The remainder of this section will make use of the following standard lemma.

\begin{lem}
\label{lemwronskian}
Let $P$ be a  polynomial with a simple zero at $a \in \C$. If the equation
\begin{equation}
P(z) y''(z) = y(z) 
\label{wronskiande}
\end{equation}
has a solution $f$ which is meromorphic in the plane and has $f(a) \in \C$, then every 
solution which  is meromorphic in the plane is a constant multiple of $f$.
\end{lem} 
\textit{Proof.} 
It may be assumed that $a=0$. The assumptions force $f(0) = 0$, and the zero of $f$ at $0$ must be simple, because otherwise $P = f/f''$ has a double zero at $a$. 
Hence 
$c_1 = f'(0) \neq 0$, 
and it follows from (\ref{wronskiande})  that $2c_2 = f''(0) \neq 0$. A second solution $g$ may then be obtained by 
writing 
$$
 \left( \frac{g}{f} \right)' (z) = \frac1{f(z)^2} = \frac1{ (c_1 z + c_2 z^2 + \ldots )^2 }
 =  \frac1{c_1^2 z^2} (1 -  2 (c_2 /c_1)  z + \ldots )  , 
 $$
and integration clearly gives rise to a logarithm. 
 \hfill$\Box$
\vspace{.1in}

\subsection{The equation (\ref{bessel})}\label{besselexample}

Let $f_4$ be as in  (\ref{besselsoln}). Then differentiating $f_4$ twice leads to 
$$
z f_4''(z) = \sum_{k=1}^\infty \frac{z^k}{(k-1)! k! } = 
\sum_{k=0}^\infty \frac{z^{k+1}}{k! (k+1)!} = f_4(z) ,
$$
after replacing $k $ by $k +1$, and so $f_4$ is a solution of (\ref{bessel}).
Lemma \ref{lemwronskian}, with $a=0$, shows that any solution of (\ref{bessel}) which is meromorphic in $\C$
is a constant multiple of $f_4$. 

It turns out that $f_4$ has a representation in terms of Bessel functions:
write $z = w^2$ and 
$$
f_4(z) =
%= \sum_{k=0}^\infty \frac{z^{k+1}}{k!(k+1)!}
 \sum_{k=0}^\infty \frac{w^{2k+2}}{k!(k+1)!} 
%= w^2 \sum_{k=0}^\infty \frac{(-1)^k (2iw)^{2k}}{2^{2k} k!(k+1)!} 
= \frac{w}{i}  \sum_{k=0}^\infty \frac{(-1)^k (2iw)^{2k+1}}{2^{2k+1} k!(k+1)!} = \frac{w}{i} J_1(2iw), $$
where $J_1$ is the Bessel function of the first kind of order $1$ \cite{Hille}. 
This relation can be used to prove that  all zeros of $f_4$ are real and non-positive,
but the following
approach applies  Green's transform \cite[pp.286-8]{Hille} directly to $f_4$ and (\ref{bessel}).

Suppose then  that 
$R > 0$ and $s \in \R$ and $R e^{is}$ is a zero of $ f_4$. Set 
\begin{equation*}
 F(r) = f_4(re^{is}) , \quad H(r) = \overline{F(r)} F'(r) .
\label{green}
\end{equation*} 
This yields, for $r > 0$, by   (\ref{bessel}),
$$
H'(r) = |F'(r)|^2 + \overline{F(r)} F''(r) =
|F'(r)|^2 + e^{2is} \overline{F(r)} f_4''(re^{is})
 =
|F'(r)|^2 + \frac{e^{is} |F(r)|^2}{r} .$$
Since $H(R) = H(0) = 0$, integration from $0$ to $R$ results in
$$
 \int_0^R |F'(r)|^2 \, dr = - e^{is}  \int_0^R \frac{|F(r)|^2}{r} \, dr,
$$
which forces $e^{is} = -1$, so that $Re^{is}$ lies on the negative real axis.

Next, a  straightforward application of the Wiman-Valiron theory \cite{Hay5} in (\ref{bessel}) shows that 
the order of $f_4$ is $1/2$, and so a standard generalisation of the Gauss-Lucas theorem
\cite{Titchmarsh} implies that all zeros of $f_4'$ are also real and non-positive. This completes the proof of the following.

\begin{lem}
\label{lembesselsoln}
The real entire function $f_4$ given by (\ref{besselsoln}) is a solution of (\ref{bessel}),
and   all zeros of $f_4$ and $f_4'$ are real
and non-positive.
%, while $f_4''/f_4$ has no zeros. 
Moreover,  any solution of (\ref{bessel}) which is meromorphic in $\C$ is a constant multiple of $f_4$. 
\end{lem}

\hfill$\Box$
\vspace{.1in}
%060323
\subsection{The equation  (\ref{hypergeometricA})} \label{hyper1a} 

Let $n \geq 2$ be an integer, and consider the equation (\ref{hypergeometricA}).
MAPLE gives solutions in terms of hypergeometric functions, 
but an explicit solution will be derived as follows.
 Let $F = F_2$ be given by (\ref{hyper1afn}) and write
 \begin{eqnarray*} 
F(z)  &=& \frac{d^{n-2} }{dz^{n-2}}  \left( z^{n-1} (z-1)^{n-1} \right) \\
&=&   \frac{d^{n-2} }{dz^{n-2}}  \left( \sum_{k=0}^{n-1} \frac{(n-1)!}{k!(n-1-k)!}  \, z^{k + n-1} (-1)^{n-1-k} \right) \\
&=& \sum_{k=0}^{n-1}  \frac{(n-1)!(k+n-1)!}{k!(n-1-k)!(k+1)!} \, z^{k + 1} (-1)^{n-1-k} \\
 &=& \sum_{k=0}^{n-1} a_k z^{k+1}  , \quad a_k \in \R, \quad a_0  \neq 0.
 \end{eqnarray*} 
 It is then clear  that $F$ is a polynomial of degree $n$, with a simple zero at $ 0$,
and that all zeros of $F$ and
$F'$ lie in $[0, 1]$, by repeated
application of the
Gauss-Lucas theorem. Moreover, $a_k$ satisfies
$$\frac{a_{k+1}}{a_k} = - \,  \frac{(k+n) (n-1-k)}{(k+1)(k+2)}   \quad \hbox{for $k =0, \ldots, n-2$.} 
$$
 Hence differentiating $F$ twice yields
\begin{eqnarray*} 
z(z-1) F''(z)  
&=&   (z^2-z) \sum_{k=0}^{n-1} (k+1)k a_k z^{k-1} \\
&=&    \sum_{k=0}^{n-1} (k+1)k a_k z^{k+1} -   \sum_{k=1}^{n-1} (k+1)k a_k z^{k} \\
&=&    \sum_{k=0}^{n-1} (k+1)k a_k z^{k+1} -   \sum_{k=0}^{n-2} (k+2)(k+1) a_{k+1} z^{k+1} \\
&=&    \sum_{k=0}^{n-1} (k+1)k a_k z^{k+1} +   \sum_{k=0}^{n-2}  (k+n)(n-1-k)  a_k z^{k+1} \\
&=&    \sum_{k=0}^{n-1} (k+1)k a_k z^{k+1} +   \sum_{k=0}^{n-1}  (k+n)(n-1-k)  a_k z^{k+1} \\
&=&    \sum_{k=0}^{n-1} (n^2-n)  a_k z^{k+1} = n(n-1) F(z) .
\end{eqnarray*}
 Thus $F$ solves (\ref{hypergeometricA}). Applying Lemma \ref{lemwronskian}, with $a=0$, 
 completes the proof of the following.
 
 \begin{lem}
 \label{lemhyper1a}
 For $2 \leq n \in \Z$, the real polynomial $F_2$ given by (\ref{hyper1afn}) has degree $n$ and solves
 (\ref{hypergeometricA}). % $$F(z) = \frac{d^{n-2}}{dz^{n-2}} \left( z^{n-1} (z-1)^{n-1} \right)$$
 Moreover, all zeros of $F_2$ and $F_2'$ are real, 
 %while $F_2''/F_2$ has no zeros, 
 and 
  any solution of (\ref{hypergeometricA}) which is meromorphic
 in the plane must be a constant multiple of $F_2$. 
 \end{lem} 
 
\hfill$\Box$
\vspace{.1in}
 
\subsection{The equation (\ref{hypergeometric3A})}\label{hyper2A}

Let $n \geq 1$ be an integer and let $K \in \R \setminus \{ 0, 1 \}$. 
%$1 < K < + \infty $ or $- \infty < K \leq -1$. 
%The following construction will lead to a rational solution of (\ref{hypergeometric3A}).
%which satisfies the hypotheses of Theorem \ref{thm1}. 
MAPLE gives solutions of (\ref{hypergeometric3A}) in terms of hypergeometric functions,
and the following direct determination of a rational
solution was found via properties of the related Jacobi polynomials \cite[p. 254]{Rainville}. 
Using the substitution 
\begin{equation}
w = \frac{K + 1}{K-1} -\frac{2K}{(K-1)z } , \quad z = \phi(w) =  \frac{2K}{K+1 - (K-1) w} ,
\quad y(z) = (z-K) h \left( w \right), 
\label{subs}
\end{equation}
write
\begin{eqnarray*}
y'(z) &=&  h\left(  w \right) + \frac{2K(z-K)}{(K-1)z^2}  h' \left( w \right),\\
y''(z) &=& \frac{2K}{(K-1) z^2}  h' \left( w \right) +
 \frac{2K}{(K-1) z^2}  h' \left( w \right) \\
 & & 
 -  \frac{4K(z-K)}{(K-1) z^3}  h' \left( w   \right) + 
  \frac{4K^2(z-K)}{(K-1)^2 z^4}  h'' \left( w  \right)\\
&=&  \frac{4K^2}{(K-1) z^3}  h' \left( w \right) +  
  \frac{4K^2(z-K)}{(K-1)^2 z^4}  h'' \left( w  \right).
\end{eqnarray*}
Observe next that, by (\ref{subs}), 
$$
 \frac{z-1}z = \frac{(K-1) (w+1)}{2K}, \quad \frac{z-K}z = \frac{(K-1)(w-1)}2.
$$
Thus substituting for $y$ and $y''$   delivers 
\begin{eqnarray*}
R(z) &=& Kn(n+1) y(z) - z^2 (z-1) (z-K) y''(z) \\
 &=& K n(n+1) (z-K) h\left(w \right)
- z^2 (z-1)(z-K) 
\left(    \frac{4K^2}{(K-1) z^3} h' \left(w \right) +  \frac{4K^2(z-K)}{(K-1)^2 z^4}  h'' \left(w \right) \right),\\
&=&  K(z-K) \left[ n(n+1)  h\left(w \right)
-  \left(  \frac{4K}{K-1} \right)  \left( \frac{z-1}z \right)  h' \left(w \right)  \right] \\
& &  -  K(z-K) \left[  \frac{4K}{(K-1)^2 } \left( \frac{(z-1)(z-K) }{z^2} \right)  h'' \left(w \right) \right] \\
&=& K(z-K) \left[   n(n+1)h(w)  - 2 (w+1) h'(w) + (1-w^2) h''(w) \right].
\end{eqnarray*}
Thus $y$ solves (\ref{hypergeometric3A}) if and only if $h$ solves 
\begin{equation}
 (1-w^2) h''(w) - 2 (w+1) h'(w) + n(n+1)h(w) = 0 .
\label{hypergeometric5A}
\end{equation}

\begin{lem}
\label{lemH}
Let
  $H(w) = H_n(w)$, with $H_n$ as  in (\ref{Hndefn}).
  Then $H$ is a polynomial of degree $n$ and solves (\ref{hypergeometric5A}). 
  Moreover, $H(-1) = 0$, all $n$ zeros of $H$ are simple, and they all lie in $[-1, 1)$. 
  \end{lem}
\textit{Proof.}
Write  (\ref{Hndefn}) in the form
\begin{eqnarray*}
H(w) 
%&=&   \frac{d^n}{dw^n} \left( (w-1)^{n-1}(w+1)^{n+1}  \right) \\
 &=&  \frac{d^n}{dw^n} \left( (w+1-2)^{n-1}(w+1)^{n+1}  \right)\\
 &=&    \frac{d^n}{dw^n} \left( \sum_{k=0}^{n-1} {n-1 \choose k} (-2)^{n-1-k} 
 (w+1)^{n+1+ k}  \right)\\ 
 &=&  \sum_{k=0}^{n-1} {n-1 \choose k} (-2)^{n-1-k} \frac{(n+1+k)!}{(k+1)!} 
 (w+1)^{ k+1}  ,
  \end{eqnarray*}  
  which leads to 
  \begin{equation}
  H(w) =  \sum_{k=0}^{n-1} b_k (w+1)^{k+1} , 
  \quad 
  b_k = 
    \frac{(n-1)! (n+1+k)!}{k!(n-k-1)! (k+1)!} (-2)^{n-1-k} ,
 \label{bkdefn0}
 \end{equation}
 in which $b_0 \neq 0$ and 
 \begin{equation}
  \frac{b_{k+1}}{b_k} = \frac{(n-k-1)(n+2+k)}{(k+1)(k+2)(-2)}
 = \frac{(k+1-n)(k+2+n)}{2(k+1)(k+2)} \quad \hbox{for $k =0, \ldots, n-2$.} 
 \label{bkrecurrence}
 \end{equation} 
  Substitution of (\ref{bkdefn0}) into the right-hand side of (\ref{hypergeometric5A}),
  followed by application of (\ref{bkrecurrence}),   delivers 
 \begin{eqnarray*}
 Q(w) &=&  (1-w^2) H''(w) - 2 (w+1) H'(w) + n(n+1)H(w) \\
 &=& (2-(w+1)) (w+1) \sum_{k=0}^{n-1} (k+1)k b_k (w+1)^{k-1} \\
 & & - 2(w+1)   \sum_{k=0}^{n-1} (k+1) b_k (w+1)^{k} + n(n+1)  \sum_{k=0}^{n-1} b_k (w+1)^{k+1}\\
 &=& 2  \sum_{k=1}^{n-1} (k+1)k b_k (w+1)^{k}   +
  \sum_{k=0}^{n-1} (n(n+1) -(k+1)k - 2(k+1))  b_k (w+1)^{k+1} \\
  &=&  2  \sum_{k=0}^{n-2} (k+2)(k+1) b_{k+1} (w+1)^{k+1}   +
  \sum_{k=0}^{n-1} (n(n+1) - (k+2) (k+1) )  b_k (w+1)^{k+1} \\
  &=&    \sum_{k=0}^{n-2} (k+1-n)(k+2+n) b_k (w+1)^{k+1}   +
  \sum_{k=0}^{n-1} (n(n+1) - (k+2) (k+1) )  b_k (w+1)^{k+1} \\   
  &=&    \sum_{k=0}^{n-1} (k+1-n)(k+2 + n) b_k (w+1)^{k+1}   +
  \sum_{k=0}^{n-1} (n(n+1) - (k+2) (k+1) )  b_k (w+1)^{k+1} \\   
  &=& 0.
  \end{eqnarray*} 
 Thus $H(w)$ is a polynomial solution of  (\ref{hypergeometric5A}),
of degree $n$, with a simple zero at $-1$, since $b_0 \neq 0$ in (\ref{bkdefn0}).
Repeated application of the Gauss-Lucas theorem to  $G(w) = (w-1)^{n-1}(w+1)^{n+1} $
shows that  all  zeros of $H(w)$ lie in 
$[-1, 1]$. Moreover, since $G$ has a zero of multiplicity $n-1$ at $1$, all zeros of 
$E = G^{(n-1)} $ lie in $[-1, 1)$ and therefore so do all zeros of $H = E'$. 
%With $D = d/dw$, Leibnitz' rule gives
%$$H(w) = \sum_{k=0}^n { n \choose k} D^k  \left( (w-1)^{n-1} \right)\timesD^{n-k}  \left( (w+1)^{n+1} \right),$$
%Since all derivatives of $(w-1)^{n-1}$, except the $(n-1)$th, vanish at $1$, while Leibnitz' rule gives
%$$
 %H(w) = \sum_{k=0}^n { n \choose k} 
%\frac{d^k}{dw^k} \left( (w-1)^{n-1} \right) \times
%\frac{d^{n-k}}{dw^{n-k}} \left( (w+1)^{n+1} \right),$$
%it follows that $H(1) \neq 0$. 
Finally, all zeros of $H$ in $(-1, 1)$ are simple,
by the existence-uniqueness theorem and  (\ref{hypergeometric5A}).
 \hfill$\Box$
\vspace{.1in}

\begin{lem} 
\label{lemF3}
With $H_n$ as in (\ref{Hndefn}) and $n \geq 1$, and with $w$ defined by (\ref{subs}),
the function $ F_3(z)  $ in (\ref{hyper2Afn}) is a rational solution 
of
(\ref{hypergeometric3A}) with $F_3(1) = F_3(K) = 0$ 
and all its zeros real and simple, and every solution of
(\ref{hypergeometric3A})  which is meromorphic in the plane is a constant multiple of $F_3$. 
Moreover, 
%$ F_3$ has the form 
$F_3(z) = P(z) z^{-n} $,
where $P$ is a real polynomial  with $P(0) \neq 0$, and $P$ has degree $n$ or $n+1$. 

Furthermore, if $P$ has degree $n+1$ and all zeros of $P$ lie in $(0 , + \infty)$, then $F_3$,
$F_3'$ and $F_3''$ have no zeros or poles in $\C \setminus \R$. In particular, this holds if $K > 1$. 
\end{lem} 
\textit{Proof.} 
 First, $F_3$ solves (\ref{hypergeometric3A}) and has a pole at $z = 0$ of order $n$,
since  $z = 0$ corresponds to $w= \infty$
and $H$ has degree $n$. 
Clearly, $F_3$ has no other poles in $\C$. 

Next, $H(-1) = 0$ and all $n$ zeros of $H$ are simple and lie in $[-1, 1)$. 
Since  one of them may be mapped to $\infty$ by $z = \phi(w)$, 
it  follows from (\ref{subs})  that $F_3(\infty) \neq 0$ and $F_3$ has $n$ or $n+1$ zeros in $\C$. 
In particular, 
$F_3$ has  zeros  at $z=1$ and  $z = K$, which correspond to $w = -1$ and $w=1$ respectively,  and
any solution of
(\ref{hypergeometric3A})  which is meromorphic in the plane is a constant multiple of $F_3$, by 
Lemma \ref{lemwronskian} with $a = K$ or $a=1$.

Now suppose that $P$ has degree $n+1$ and all zeros of $P$ lie in $(0 , + \infty)$.  
This will certainly hold if $K > 1$, 
because in this case
the function 
$z = \phi(w)$ is finite and  increasing for $-1 \leq w < 1$,  and maps 
$[-1, 1)$ to $[1, K)$, so that 
$F_3(z) = (z-K) H\left(w \right)$ inherits all $n$  zeros of $H(w)$, as well as having a zero
at $z=K$. 
Under these assumptions,  a consideration of leading terms shows that $zP'(z) - n P(z) $ has degree $n+1$, and so
$$
F_3'(z) = \frac{zP'(z)- nP(z)}{z^{n+1}}
$$
has $n+1$ zeros in $\C$. Of these, $n$ arise from Rolle's theorem and lie in $(0, + \infty )$, 
while  one lies in $(-\infty, 0)$ because, with $x \in \R$, 
$$
\frac{F_3'(x)}{F_3(x) }  \sim \frac1x < 0 \quad \hbox{as $x \to - \infty$,} \quad
\frac{F_3'(x)}{F_3(x) }  \sim - \, \frac{n}x > 0 \quad \hbox{as $x \to 0-$.}$$
Thus $F_3$ and $F_3'$ have no zeros in $\C \setminus \R$, and nor has $F_3''$, because of 
 (\ref{hypergeometric3A}). 
\hfill$\Box$
\vspace{.1in}
%060323

Taking $K=2$ delivers $F_3(z) = (z-K)H_n(w) = (z-2)H_n(3-4/z)$ and (with help from MAPLE) the following:
\begin{eqnarray*}
&n=1,& \quad g_1(z) =  \frac{8(z-1)(z-2)}z;\\
&n=2,& \quad g_2(z) = \frac{144(z-1)(z-4/3)(z-2)}{z^2};\\
&n=3,& \quad g_3(z) = \frac{384(z-1)(z-2)( 11z^2-30z+20 )}{z^3} .
\end{eqnarray*}
In these examples, $g_j$, $g_j'$  and $g_j''$ have no zeros in $\C \setminus \R$, by Lemma \ref{lemF3}.
On the other hand, 
choosing $K = -1$ leads to  
$ F_3(z) = (z-K)H_n(w) = (z+1)H_n(-1/z)$, as well as: 
\begin{eqnarray*}
&n=1,& \quad h_1(z) = \frac{2(z^2-1)}z;\\
&n=2,& \quad h_2(z) = \frac{-12(z^2-1)}{z^2};\\
&n=3,& \quad h_3(z) = \frac{-24( z^2 - 1)( z^2-5)}{z^3} ;\\
&n=4,& \quad h_4(z) = \frac{720(z^2-1)( z^2-7/3)}{z^4} .
\end{eqnarray*}
Here $h_j$, $h_j'$ and $h_j''$ have  no zeros in $\C \setminus \R$
for $j=2, 4$, but $h_j'$ has non-real zeros for $j=1, 3$. 
%The difference 
%between the cases of odd and even $n$ reflects the fact that $H_2$ and $H_4$ vanish at $w=0$, which corresponds
%to $z = \infty$, but $H_1$ and $H_3$ do not. 

\subsection{The equation (\ref{hypergeometric3AB})}\label{hyper2AA}

A solution of (\ref{hypergeometric3AB}) is obtained by the following limiting process with $K$ real: let
$n \geq 1$ and let $F_3, H_n$ be as in (\ref{hyper2Afn}) and (\ref{Hndefn}), and set
$$
F_4(z) = \lim_{K \to + \infty} \frac{F_3(z)}{-K} =  \lim_{K \to + \infty}
\left( \frac{z-K}{-K} \right)  H_n \left(  \frac{K + 1}{K-1} -\frac{2K}{(K-1)z } \right) =
H_n \left( 1 - \frac2z \right). 
$$
Since all zeros of $H_n$ and $H_n'$ lie in $[-1, 1)$, by Lemma \ref{lemH},  $F_4$ and $F_4'$ have no zeros in
$\C \setminus [1, + \infty)$, 
and $F_4$ has a pole of order $n$ at $0$. Applying Weierstrass' theorem
% in a neighbourhood of $z = i$ 
yields,
since $F_3$ solves    (\ref{hypergeometric3A}),
$$
\frac{F_4''(z)}{F_4(z)} = \lim_{K \to + \infty} \frac{F_3''(z)}{F_3(z)} = 
\frac{-n(n+1)}{z^2(z-1)} 
$$
as required. Furthermore, $F_4$ has a simple zero at $z = 1$, inherited from the simple zero of $H_n$ at $w = -1$, 
which completes the proof of the following. 

\begin{lem} 
\label{lemF4}
With $H_n$ as in (\ref{Hndefn}) and $n \geq 1$, the function $ F_4(z) =  H_n \left(1 - 2/z \right) $ is a rational solution 
of
(\ref{hypergeometric3AB}) with a simple zero at $1$, and every solution of
(\ref{hypergeometric3AB})  which is meromorphic in the plane is a constant multiple of $F_4$. 
Furthermore,  $F_4$,
$F_4'$ and $F_4''$ have no zeros or poles in $\C \setminus \R$. 
\end{lem} 
\hfill$\Box$
\vspace{.1in}

Calculating $F_4(z) = H_n(1-2/z)$ using MAPLE delivers: 
\begin{eqnarray*}
&n=1,& \quad p_1(z) =  \frac{4(z-1)}z;\\
&n=2,& \quad p_2(z) = \frac{24(z-1)(z-2)}{z^2};\\
&n=3,& \quad p_3(z) = \frac{192(z-1)( z^2 -5z+5)}{z^3} .
\end{eqnarray*}
%It is easy to verify that $p_j$, $p_j'$  and $p_j''$ have only real zeros for $j=1, 2, 3$. 

\section{Lemmas needed for the proof of Theorem \ref{thm1}} 

\begin{lem}
\label{lemrealzeros}
Let $g$ be a real meromorphic  function in the plane, of order at most $1$, and with infinitely many zeros, all but finitely many of them real, and assume that $g$ has finitely many poles.
Then 
$$
\lim_{y \to +\infty, y \in \R} 
\frac{\log |g(iy)|}{\log y} = + \infty . 
$$
\end{lem} 
\textit{Proof.} It is enough to prove this when $g$ is real entire, with only real zeros,  and with $g(0) \neq 0$.
The hypotheses then imply that 
$$
g(z) = e^{\alpha z + \beta } \prod_{n=1}^\infty \left( 1 - \frac{z}{a_n} \right) e^{z/a_n} ,
$$
with $\alpha, \beta , a_n$ real. As $y \to +\infty $ with $ y \in \R$ this gives
\begin{eqnarray*}
2 \log |g(iy)| &\geq& 
2 \sum_{n=1}^\infty \log \left| 1 - \frac{iy}{a_n} \right|  - O(1) 
=
 \sum_{n=1}^\infty \log \left( 1 + \frac{y^2}{a_n^2} \right)  - O(1) \\
&\geq & 
 \sum_{|a_n| \leq \sqrt{y} } \log \left( 1 + y \right)  - O(1) 
\geq n(\sqrt{y} , 1/g) \log y - O(1) . 
\end{eqnarray*}
\hfill$\Box$
\vspace{.1in}

\begin{lem}[\cite{BLa}]\label{lemnorfam1}
Let $D \subseteq \C$ be a domain and let $\mathcal{F}$ be a family of meromorphic functions $f$ on $D$ such that $f$ and $f''$ have no zeros in $D$. Then the family $\{ f'/f : \, f \in \mathcal{F} \}$ is normal on $D$.
\end{lem}
\hfill$\Box$
\vspace{0.1in}

Next,  suppose that  $G$ is a transcendental meromorphic 
function in the plane,
and that  $G(z) \to a \in \C \cup \{ \infty \}$ as $z \to \infty$ along a path
$\gamma $; then the inverse $G^{-1}$ is
said to have a transcendental singularity
over the asymptotic value $a$~\cite{BE,Nev}. If $a \in \C$ then for each $\varepsilon > 0$ 
there  exists a component $\Omega = \Omega( a, \varepsilon, G)$ of the
set $\{ z \in \C : |G(z) - a | < \varepsilon \}$ such that
$\gamma \setminus \Omega$ is bounded: these components are referred to as neighbourhoods of the singularity \cite{BE}.
Two  such paths $\gamma, \gamma'$ on which $G(z) \to a$ determine distinct singularities if the corresponding components
$\Omega( a, \varepsilon, G)$, $\Omega'( a, \varepsilon, G)$ are disjoint for  some $ \varepsilon > 0$.
The singularity 
%of $G^{-1}$ corresponding to $\gamma$
is called direct  \cite{BE} if $\Omega( a, \varepsilon, G)$, for some
$\varepsilon > 0$, contains finitely many zeros of $G-a$, and 
%In the contrary case the singularity is called direct and 
indirect otherwise. A direct singularity is called logarithmic if there exists $\varepsilon > 0$ 
such that $w = \log 1/(G(z)-a)$ is a conformal bijection from $\Omega( a, \varepsilon, G)$ to the half-plane
${\rm Re} \, w > \log 1/\varepsilon $.
%, in which case $\Omega( a, \varepsilon, G)$, for all sufficiently small $\varepsilon > 0$, contains no zeros of $G-a$. 
Finally, transcendental singularities over $\infty$ may be 
classified using $1/G$, and a transcendental singularity will be referred to as lying in an open set $D$ 
%upper half-plane $H = \{ z \in \C : \mathrm{Im} \, z >  0 \}$ 
if $\Omega(a, \varepsilon, G) \subseteq D$ for some $\varepsilon > 0$.

The next lemma combines [\cite{Lajda09}, Lemma 2.4] and  [\cite{Lawiman13}, Lemma 2.2], and 
links asymptotic values approached on paths in the upper half-plane $H^+$  with the growth of 
the Tsuji characteristic $ \mathfrak{T} (r, g) =  \mathfrak{m} (r, g)+  \mathfrak{N} (r, g)$  for functions $g$ that are  meromorphic
on the closed upper half-plane \cite{BEL,GO,Tsuji}.
%$\overline{H} =\{ z \in \C : \mathrm{Im} \, z \geq  0 \}$.

\begin{lem}
[\cite{Lajda09,Lawiman13}]
\label{directsinglem2}
Let $L \not \equiv 0$ be a real  meromorphic function in the plane such that
$\mathfrak{T}(r, L) = O( \log r )$ as $r \to \infty$, and let $F(z) = z-1/L(z)$. 
Assume that at least one of $L$ and $1/L$ has finitely many non-real poles. 
Then there exist finitely many $\alpha \in \C $ such that $F(z)$ or $L(z)$
tends to $\alpha$ as $z$ tends to infinity along a path in $\C \setminus \R$.

Moreover,  there exists 
at most one direct transcendental singularity of $F^{-1}$ lying in 
$H^+$.
\end{lem}
\hfill$\Box$
\vspace{0.1in}

The following result of Levin and Ostrovskii \cite{LeO} 
(see also \cite[ Ch. 6, Lemma 5.2]{GO} and  \cite[Lemma 2.4]{Lajda09}) will be required. 

\begin{lem}[\cite{LeO}]\label{lemlevinost}
Let $G$ be a meromorphic function in the plane: then, for each $R \geq 1$,
$$
\frac1{2 \pi} \int_R^{+ \infty} \frac1{r^3} \int_0^\pi \log^+ |G(re^{i \theta })| \, d \theta \, dr 
\leq 
\int_R^{ + \infty } \frac{  \mathfrak{m} (r, G)}{r^2} \, dr .$$
If, in addition, $G$ is real meromorphic with finitely many poles, and satisfies 
$\mathfrak{T}(r, G) = O( \log r )$ as $r \to \infty$, then $T(R, G) = O( R \log R )$ as $R \to + \infty$. 
\end{lem}
\hfill$\Box$
\vspace{0.1in}

\begin{lem}
\label{lemlevin}
There exists a positive constant $c_0$
such that if the function $\psi $ maps the upper half-plane $H^+$ analytically into itself then, for 
 $ r \geq 1$ and 
$\theta\in (0,\pi)$,
\begin{equation}
\frac{|\psi(i)| \sin\theta}{5 r} < |\psi(re^{i\theta})| <
\frac{5r |\psi(i)|}{\sin\theta}\quad \hbox{and}  \quad 
 \left| \frac{ \psi'( re^{i \theta } ) }{ \psi( re^{i \theta } ) } \right|
 \leq \frac{c_0}{r \sin \theta} .
  \quad\label{C1}
\end{equation}
\end{lem}
Both of these estimates  are standard: the first is essentially just Schwarz' lemma
 \cite[Ch. I.6, Thm $8'$]{Le}, while the second follows from Bloch's theorem applied to $\log \psi$.
\hfill$\Box$
\vspace{0.1in}

\subsection{The Levin-Ostrovskii factorisation}\label{levostfact} 

The following constructions are standard \cite{BEL, LeO}. Suppose that $(u_k), (v_k)$ are sequences 
%defined and 
satisfying 
$u_k < v_k < u_{k+1}$
for $- \infty \leq M < k < N \leq + \infty$. Then there exists $k_0 \in \N$ such that 
 $u_k$ and $ v_k$ have the same sign for $|k| \geq k_0$, and
  $$
  \psi (z) =  \prod_{|k| \geq k_0}  \frac{1-z/v_k}{1-z/u_k} 
 $$
 converges on $\C$ by
 the alternating series test. 
 Furthermore, $\psi$ satisfies,  for $z$ in the upper half-plane $H^+$,
$$
\arg \psi(z) =
\sum_{|k| \geq k_0}  \arg \frac{1-z/v_k}{1-z/u_k} = 
\sum_{|k| \geq k_0}  \arg \frac{v_k-z}{u_k-z} \in (0, \pi ) .
$$

This leads to  the Levin-Ostrovskii factorisation \cite{BEL,LeO} of the logarithmic derivative
of a real entire  function $f$ with real zeros. If $f$ has finitely many zeros,
set $\psi (z) = 1$, while if $f$ has infinitely many zeros $u_k$ then zeros of $f'$ given by Rolle's theorem can be labelled $v_k$ so that $u_k < v_k < u_{k+1}$, whereupon  $\psi$ may be constructed as above. It follows that $f'/f = P \psi $, where $P$ is  real meromorphic  with finitely many poles and either 
$\psi \equiv 1$ or $\psi (H^+) \subseteq H^+$.

\section{Proof of Theorem \ref{thm1}: first steps} 

Let $f$ be as in the hypotheses and write
\begin{equation}
\label{ww1}
L = \frac{f'}f, \quad F(z) = z - \frac{f(z)}{f'(z)}, \quad F' = \frac{ff''}{(f')^2}  .
\end{equation}

\begin{lem}
\label{lemnorfam}
Let $0 < \delta  < \pi /2$ and $\delta < \sigma < \pi - \delta $.\\
(I) If $rL(re^{i \sigma } )$ is bounded as $r \to + \infty$ then $zL(z)$ is  bounded  
as $z \to \infty $ with $\delta < \arg z < \pi - \delta $.\\
(II)  If $\lim_{r \to +\infty} rL(re^{i \sigma } ) = 0$, then 
$zL(z) \to 0$ uniformly as $z \to \infty $ with $\delta < \arg z < \pi - \delta $.
\end{lem}
\textit{Proof.}
The functions $u_R(z) =
R L(Rz)$, $R \geq 1$, form a normal family on the domain
$D_1 = \{ z \in \C : \, 1/2 < |z| < 2, \, \delta /2 < \arg z < \pi - \delta /2 \}$: this follows from
Lemma \ref{lemnorfam1} applied to the functions $f(Rz)$.
Take a sequence $R_n \to + \infty$ such that
$(u_{R_n} )$ converges locally spherically uniformly on $D_1$. 
In case (I),  $(u_{R_n} )$ cannot have $\infty$ as limit,
while  
in case (II) the limit function must vanish identically, 
by the identity theorem. 
\hfill$\Box$
\vspace{.1in}

\begin{lem}
\label{lemww1}
Poles of $F$ in $\C$ coincide with  zeros of $L=f'/f$, all of which are real and simple. 
All zeros of $F' $ in $\C$ 
are real  zeros of $f$ and super-attracting fixpoints of $F$; furthermore, simple zeros of $F'$ in $\C$ are zeros 
of $f$ which are not zeros of $f''$, while multiple zeros of $F'$ in $\C$ have multiplicity $2$ and are common simple zeros
of $f$ and $f''$. 
\end{lem}
\textit{Proof.}
This is standard, and  all assertions follow from (\ref{ww1}). First, any
multiple zero of $L=f'/f$ would be a zero of $f''$, and hence of $f$, and thus a pole of $f'/f$, an obvious contradiction.
Next, zeros of $F'$ are zeros of $f$ or $f''$, and hence of $f$. But multiple zeros of $f$ are not zeros of $F'$, and so
all zeros of $F'$ must be simple zeros of $f$, and since $f''/f$ has no zeros they cannot be 
zeros of $f''$ of multiplicity greater than $1$. 
% since zeros of $f''$ are zeros of $f$ and hence fixpoints of $F$, while% and hence of $f$. 
\hfill$\Box$
\vspace{.1in}

Define the sets  $W^+$ and $W^-$ using 
\begin{equation}
H^+ = \{ z \in \C : \,   {\rm Im} \, z > 0 \}, \quad 
%H^- = \{ z \in \C : \,   {\rm Im} \, z < 0 \}, \quad 
W^{\pm} = \{ z \in H^+: \, \pm F(z) \in H^+ \}.
\label{ww2}
\end{equation}

The next lemma is fairly standard and goes back to 
Sheil-Small \cite{SS}. 

\begin{lem}
\label{lemww7}
Let $x_0 \in \R$ be a zero of $f'/f$. If  $(f'/f)' (x_0) < 0$
then  $x_0  \in \partial W^-$, while if $(f'/f)' (x_0) > 0$
then  $x_0  \in \partial W^+$.
Poles of $f$ are repelling  fixpoints of $F$ and lie in $\partial W^+ \setminus \partial W^-$. 
\end{lem}
\textit{Proof.} The first two assertions hold since as $z \to x_0$ from within $H^+$ the sign 
of ${\rm Im} \, (- f(z)/f'(z))$ is the same as that of 
${\rm Im} \, ( f'(z)/f(z))$. Furthermore, if $x_1$ is a pole of $f$ of multiplicity $m_1$ then 
$F(x_1) = x_1$ and $F'(x_1) = 1 + 1/m_1 > 1$. 
%here exists $c \in \R$ with 
%$f'(z)/f(z) \sim - c^2 (z-x_0) $  and $F(z) \sim -f(z)/f'(z) \sim 1/c^2 (z-x_0) $ as $z \to x_0$.  
\hfill$\Box$
\vspace{.1in}

\begin{lem}
The following statements hold.\\
(i) If $F$ is transcendental and has finitely many asymptotic values then 
all but finitely many zeros of $f$  are simple.\\
(ii) 
If $F$ is rational and either $F(\infty) = \infty $ or $\infty$ is not a multiple point of $F$, 
%$\infty$ is either a simple zero of $F$ or a multiple pole of $F$, 
then all zeros of $f$ in $\C$ are simple. 
%(ii) All zeros of $L$ are simple, as are all poles of $F$. \\
%(iii) All but finitely many zeros of $F'$ have multiplicity $1$ or $2$.
\label{lemww6}
\end{lem} 
\textit{Proof.} 
This uses standard facts involving iteration \cite{Stei}. 
To prove (i) observe that 
a multiple zero of $f$ is an attracting, but not super-attracting, fixpoint of $F$, and so under iteration of $F$
attracts a critical or asymptotic value of $F$, while zeros of $F'$ in $\C$ are fixpoints of $F$. 
Now (ii) follows since the only singular values of $F^{-1}$ are the values taken by $F$ at multiple points in
$\C \cup \{ \infty \}$, all of which are fixpoints of $F$ by Lemma \ref{lemww1} and the assumptions of (ii). 
 %Next, a zero of $L$ is a zero of $f'$ which is not a zero of $f$, and hence not a zero of $f''$,
 %and is therefore simple. 
%Assume for simplicity that $f$ has only simple zeros: then so has $f'$, because zeros of $f''$ are zeros of $f$, and so has $f'/f$. 
%Finally, (iii) holds because all but finitely many zeros of $F'$ are simple zeros of $f$. 
\hfill$\Box$
\vspace{.1in}

Denote by $\partial D$ the boundary  of a domain $D$ with respect to $\C$.
\begin{lem}
\label{lemww8a1}
Let $C, D$ be domains with $C \subseteq D \subseteq H^+$ and $\R \subseteq \partial D$,
such that $F$ maps 
$C$ univalently onto $D$. Then  $\partial C$ 
contains at most one point which is a pole of $f$. 
\end{lem}
\textit{Proof.} Suppose that $y_1, y_2 \in \partial C$ are distinct poles of $f$. Each $y_j$ is a real
repelling fixpoint of $F$ and the branch of $F^{-1}$ mapping $D$ to $C$  extends
by reflection  to a small 
neighbourhood $U_j$ of $y_j$, with an attracting fixpoint at $y_j$. The iterates $(F^{-1})^n$ 
of $F^{-1}: D \to C \subseteq D \subseteq H^+$  form 
a normal family on $D$, but as $n \to \infty$ they tend
to the constant $y_j$ on $D \cap U_j$, a contradiction.
%which is impossible. 
\hfill$\Box$
\vspace{.1in}

\begin{lem}
\label{lemww8a}
Let $A$ be a  component of $W^+$,
and suppose that a closed interval $[a, b]$ lies in  $\partial A \cap \R$, with $a < b$ 
and $f(a),  f(b) \in \{ 0 , \infty \}$, and with 
$f(x) \neq 0, \infty$ on $(a, b)$.  Then one of the following holds: \\
(A) $f(a) \neq f(b)$ and $L = f'/f$ has no zeros in $(a, b)$; \\
(B) $f(a) = f(b) = \infty$, the function $L$ has exactly one zero $c$ in $(a, b)$, and $c$ satisfies
$L'(c) > 0$, while $F$ does not map $A$ univalently onto  $H^+$.  
\end{lem}
\textit{Proof.} 
Observe first that all zeros of $L$ in $\C$ are simple, by Lemma \ref{lemww1}, 
and that  if $f(a) = f(b) = \infty$ then $F$ cannot map $A$ univalently onto $H^+$, by 
Lemma \ref{lemww8a1}. Moreover, if $f(a) \neq f(b)$ then $L = f'/f$ has an even number of zeros in $(a, b)$.
It follows that if neither (A) nor (B) holds then there exists
at least one 
zero $d$ of $L$ in $(a, b)$ with 
$L'(d) <  0$, contradicting Lemma \ref{lemww7}.  

\hfill$\Box$
\vspace{.1in}

\begin{lem}
\label{lemww8}
Let $A$ be a  component  of $W^+$ which is mapped univalently onto $H^+$ by $F$, and assume that 
$x_1 \in \partial A \cap \R$ is a zero of $L = f'/f$. 
Then at least one of $(-\infty, x_1]$ and $[x_1, + \infty)$ lies in $\partial A$. 
\end{lem}
\textit{Proof.}  Assume the contrary. Since all multiple points of $F$ in $\C$
are zeros of $f$, by Lemma \ref{lemww1},
 it is possible to start at $x_1$ and follow $\R$ in each direction until the first encounter with a zero or pole of $f$, giving a closed interval 
$[a, b] \subseteq \partial A \cap \R$, with $a < x_1 < b$, satisfying the hypotheses
of Lemma \ref{lemww8a}; this is impossible, since 
 alternative (A) is  incompatible with the existence of
$x_1$ and (B) with $F$ mapping $A$ univalently onto $H^+$. 
\hfill$\Box$
\vspace{.1in}

\begin{lem}
\label{lemww8aa}
Let $A$ be a  component  of $W^+$. Then $A$ is unbounded.
\end{lem}
\textit{Proof.}  Assume the contrary: since $F$ has no critical values in $\C \setminus \R$,
the mapping $F: A \to H^+$ is univalent and onto. Thus $F$
must have a pole on $\partial A \cap \R$, which contradicts Lemma~\ref{lemww8}.
\hfill$\Box$
\vspace{.1in}

\begin{lem}
\label{lemww8aaa}
Let $A$ be a  bounded component  of $W^-$. Then $-F$ maps $A$ univalently onto $H^+$,
and $\partial A$ 
 consists of a closed interval $[a, b]$, where 
$-\infty < a < b < + \infty$ and $f(a) = f(b) = 0$, 
 together with a Jordan curve 
 $\lambda$ which joins $a$ to $b$ via  $H^+$. Moreover, 
  $\partial A$ contains precisely one zero
 $x_0 \in (a, b)$ of $L = f'/f$.
 \end{lem}
\textit{Proof.}  
First, $F$ must have a pole on $\partial A$, and so on $\partial A \cap \R$. Second, the mapping is univalent since $F$ has no critical values in $\C \setminus \R$. Finally, the nature of the boundary follows
from the absence of  bounded components of $W^+$. 
\hfill$\Box$
\vspace{.1in}

%070323

\begin{definition}
\label{defnchain} 
A finite chain $D$ of bounded components of $W^-$ will mean the following:\\
(a) $D$ is the union of $N \in \N$  bounded components $C_{1}, \ldots, C_{N}$ of $W^-$, each as in Lemma \ref{lemww8aaa};
 \\
(b) the boundary of each $C_{j}$ consists of a closed interval $[a_{j}, b_{j}]$, where 
$-\infty < a_j < b_j < + \infty$, 
 together with a Jordan curve 
 $\lambda_{j} $ which joins $a_{j}$ to $b_{j}$ via  $H^+$;\\
 (c) 
 the boundaries of the $C_{j}$ are disjoint except  that 
 $b_{j-1} = a_{j}$. 
 
 Such a finite chain $D$ 
 will be called maximal if $D' = D$ whenever $D'$ is a finite chain of bounded components of $W^-$  with $D \subseteq D'$. 
 \end{definition}
   
  \begin{lem}
  \label{lemchain}
  Let $D$ be a maximal finite chain  of bounded components of $W^-$ as in Definition~\ref{defnchain}. 
  Then   $a_2 = b_1, \ldots , a_{N} = b_{N-1}$ are common simple zeros 
  of $f$ and $f''$, and double zeros of $F'$, and  there exists a component $A$ of $W^+$  such that 
  $$
    \lambda_1 \cup \ldots \cup \lambda_N \subseteq \partial A.
  $$
  Moreover, 
 if $x^* = a_1$ or $x^* = b_N$ then $x^*$ satisfies exactly one of the following:  (i)
  $x^*$ is  a simple zero of $F'$ and a simple zero of $f$, but not a zero of $f''$; (ii)
  $x^*$ is a double zero of $F'$,  and a common zero of $f$ and $f''$,
  lying on the boundary of an unbounded component $B$ of $W^-$.
  %Moreover, there exist a component $A$ of $W^+$ and $c_1 < a_1, d_N > b_N$, such that 
 % $$[c_1, a_1] \cup  \lambda_1 \cup \ldots \cup \lambda_N \cup [b_N, d_N] \subseteq \partial A.$$
  \end{lem}
\textit{Proof.}  This follows from the maximality of $D$ and Lemmas \ref{lemww1} and \ref{lemww8aaa}. 
\hfill$\Box$
\vspace{.1in}

\section{The case where $f$ is a rational function}\label{rationalcase} 

\begin{prop}
Assume that $f$ is a rational function which satisfies 
the hypotheses of Theorem~\ref{thm1}. Then 
 %either $f$ is a polynomial of degree $2$ or 
 there 
 exist $\alpha_1 , \alpha_2 , \alpha_3 \in \R$ 
 with $\alpha_1 \alpha_2 \neq 0$
 such that $g(z) = \alpha_1 f(\alpha_2 z+ \alpha_3)$ is one of the functions $F_j$ in (v), (vi), (vii) and (viii) 
 of Theorem \ref{thm1}.
% \begin{eqnarray*}
% &(I)& \quad z^Q, \quad \hbox{ for some $Q \in \Z$} ; \\ &(II)& \quad 
 %\frac{d^{n-2} }{dz^{n-2}} \left( z^{n-1} (z-1)^{n-1} \right), \quad  \hbox{for some $n \geq 2$};\\&(III)& \quad 
%(z-K) H \left(  \frac{K + 1}{K-1} -\frac{2K}{(K-1)z } \right) \quad 
 %\hbox{for some $K \geq 1$, where $H$ is given by (\ref{Hdefn})};\\&(IV)& \quad a + z^{-Q} ,\quad 
 %\hbox{for some $a \in \R$ and  $Q = 1, 2$}.\end{eqnarray*} 
 \label{proprational}
 \end{prop} 
 
 The whole of this section will be occupied with the proof of Proposition \ref{proprational}.
 First, $f/f''$ is a polynomial and 
\begin{equation}
\label{wwrat1}
L(z) = \frac{f'(z)}{f(z)} = \frac{m}z + O( |z|^{-2} )  \quad \hbox{as $z \to \infty$, }
\end{equation}
where $m $ is the number of zeros minus the number of poles, counting multiplicities, of $f$ in the finite plane. Further, $F(\infty)$ exists and is real or infinite, and all components of $W^\pm $ are mapped univalently onto $H^+$ by $\pm F$.

\subsection{The case where $m \neq 0, 1$} 

Assume that $m \neq 0, 1$ in (\ref{wwrat1}): then
\begin{equation}
\frac{f''(z)}{f(z)} = L'(z) + L(z)^2 = \frac{m(m-1)}{z^2} + O( |z|^{-3} ) \quad \hbox{ as $z \to \infty$. } 
\label{mnot01}
\end{equation}
Thus 
$f/f''$ has degree $2$ and so has either
one real double zero, or two simple real zeros.

%\subsubsection{The case where $m \neq 0, 1$ and $f/f''$ has a double zero}

Suppose that $f/f''$ has a double zero. Then applying a real translation in the $z$ plane leads to
 $f(z)/f''(z) = c z^2 $ for some real constant $c$, and comparison 
with (\ref{mnot01}) forces $f$ to satisfy
$$
z^2 y''(z) = m(m-1)  y(z),
$$
which has linearly independent solutions $z^{d_j}$, where $d_1 = m \neq 0, 1$ and $d_2 = 1-m \neq 0, 1$. 
If $f$ is a constant multiple of $z^{d_j}$, for some $j$, then clearly  $f$ satisfies
conclusion (v) of Theorem \ref{thm1}. 
The only remaining possibility in this subcase  is that  there exist $a_1, a_2 \in \C \setminus \{ 0 \}$ with
$$
f(z) = a_1 z^{m} + a_2 z^{1-m} = z^m ( a_1 + a_2 z^{1-2m} ). 
$$
Since $f$ has only real zeros,  the odd integer  $1-2m$ 
must be  $ \pm 1$,
and either possibility  gives $m=0$ or  $m=1$, a contradiction. 

%\hfill$\Box$\vspace{.1in}\subsubsection{The case where $m \neq 0, 1$ and $f/f''$ has two simple zeros}

Assume next that $f/f''$ has two simple zeros, which implies that  $f$ has no poles in $\C$ and 
%one possibility is clearly that $f$ is a polynomial of degree $2$ with real zeros. 
%Assuming this  not to be the case,  
is a polynomial of degree $n \geq 2$.
A real linear change of variables then leads to 
$$
z (z-1)  f''(z) = d f(z), \quad d \in \R \setminus \{ 0 \} .
$$
A  comparison of leading terms shows that $d = n(n-1)$, giving equation (\ref{hypergeometricA}),
so that $f$ satisfies 
 conclusion (vi) 
of Theorem \ref{thm1}, by Lemma \ref{lemhyper1a}.

\hfill$\Box$
\vspace{.1in}

\subsection{The case $m = 1$}

Suppose that $m = 1$ in (\ref{wwrat1}): if $f$ has no poles in $\C$ then evidently $f$ is a linear function, contradicting the assumption that $f''/f$ has no zeros. Assume for the remainder of this
section that $f$ has at least one pole in $\C$. Then a real linear re-scaling delivers $c \in \R \setminus \{ 0 \} $ and $q \geq 1$
such that, as $z \to \infty$ and $\zeta = 1/z \to 0$, 
$$
f(z) =  z \left( 1 + \frac{c}{z^{q+1} } + \ldots \right) ,
\quad
J(\zeta ) = \frac1{f(1/\zeta )} = \zeta  ( 1 - c \zeta^{q+1} + \ldots ) = \zeta - c \zeta^{q+2} + \ldots .
$$
The flower theorem from complex dynamics  \cite{Stei} (see  \cite[Lemma 10]{BEL2}  for a convenient statement of  the theorem as applied here)
gives $q+1$ components $U_j$ of the Fatou set of $J$, each with $0 \in \partial U_j$ 
and containing a critical value $\zeta_j$ of $J$, such that the iterates $J^n$ tend to $0$ on $U_j$. Moreover  the $U_j$
can be labelled so that, as $n \to + \infty$, 
$$
\arg J^n(\zeta_j) \to \frac{2 \pi j - \arg c }{q+1} .
$$
Since all critical values of $J$ belong to $\R \cup \{ \infty \}$, because those of $f$ do, while
$J( \R \cup \{ \infty \} ) \subseteq \R \cup \{ \infty \}$, this gives  a contradiction
unless $q=1$.

It follows that, again as $z \to \infty$,
\begin{eqnarray}
f(z) &=& z + c/z + \ldots, \quad 
f'(z) = 1 - c/z^2 + \ldots, \quad f''(z) = 2c/z^3 + \ldots , \nonumber \\
\frac{f(z)}{f'(z)} &=& \frac{z( 1+c/z^2 + \ldots )}{ 1 - c/z^2 + \ldots } 
= z ( 1 + 2c/z^2 + \ldots ) = z + 2c/z + \ldots , \nonumber \\
 F(z) &=& -2c/z + \ldots , \quad \frac{f''(z)}{f(z)} = 2c/z^4 + \ldots .
 \label{Finfty}
\end{eqnarray} 
Hence $W^+$ and $W^-$ have one 
unbounded component $A$ between them. 
Moreover,  all zeros of $f$ are simple by Lemma \ref{lemww6}. 

Since $f$ has at least one pole in $\C$, 
Lemmas \ref{lemww7}, \ref{lemww8a1} and \ref{lemww8aa} imply that $A$ is a component of $W^+$ and  $f$ has exactly 
 one pole $x_0$ in $\C$, of order $n$ say, and hence $n+1$  zeros in $\C$, 
 all simple. 
 Furthermore, each of these
 simple zeros $u$ of $f$  
 %either is a simple zero of $f/f''$ or has $f(u)/f''(u) \neq 0, \infty$, and each $u$ 
 is a multiple point of $F$ and so lies in $\partial W^+ \cap \partial W^-$. Because $W^-$ has only bounded components, 
 each $u$ belongs to the boundary of a  maximal finite chain $D$ of bounded
components of $W^-$ as in Definition \ref{defnchain}.

Take such a maximal finite chain $D$: then the unique pole $x_0$ of $f$ does not lie on $\partial D$, by Lemma~\ref{lemww7}.
Moreover, with the notation of  Definition \ref{defnchain} and Lemma \ref{lemchain}, $a_1$ and $b_N$ are simple poles of $f''/f$, but the intermediate points $a_j = b_{j-1}$, $j= 2, \ldots N$, are neither zeros nor poles of $f''/f$. Hence the closure of
each such $D$ contributes exactly $2$ to the number of poles of $f''/f$ in $\C$. Since $f''/f$ has a double pole at $x_0$, and no zeros in $\C$, (\ref{Finfty}) implies that there exists precisely one maximal finite chain~$D$.

Hence among the $n+1$ zeros of $f$, precisely
$n-1$ are also zeros of $f''$, and   $f/f''$ has degree $4$ and two simple zeros at the ends $a_1, b_N$ of $D$,
plus  one double zero at $x_0$. Since $x_0 \in \partial A \setminus \partial D$, a real linear change of independent variable 
makes it possible to assume that $x_0 = 0$ and $\partial D \cap \R = [1, K]$ for some $K > 1$. 
Hence $f$ satisfies
$$
z^2 (z-1) ( z-K) f''(z) = d f(z) , \quad d \in \R ,
$$
and expanding about $z=0$ shows that $d = Kn(n+1)$, giving equation (\ref{hypergeometric3A}).  Lemma \ref{lemF3} then implies that $f$ satisfies conclusion (vii) of Theorem \ref{thm1}. This completes the discussion of the case $m=1$. 
\hfill$\Box$
\vspace{.1in} 
 
\subsection{The case $m=0$} 

Suppose that $m=0$ in (\ref{wwrat1}): then  $f$ has as many zeros as poles in $\C$, counting multiplicities. Moreover, 
$f(\infty) $ is finite and real but non-zero, and  it may be assumed that $f(\infty) = 1$. 
Further,  
there exist $c \in \R \setminus \{ 0 \} $ and $s \geq 1$ such that, as $z \to \infty$,
\begin{equation}
f(z) - 1 \sim c z^{-s} , \quad 
L(z) = \frac{f'(z)}{f(z)} \sim f'(z) \sim -cs z^{-1-s} , \quad F(z) \sim \frac{z^{1+s}}{cs} , 
\label{wwrat3a}
\end{equation}
as well as
\begin{equation}
 \frac{f''(z)}{f(z)}
= L'(z) + L(z)^2  \sim cs(s+1) z^{-2-s} .
\label{wwrat4a}
\end{equation}
Thus $F$ has a pole of multiplicity $1+s \geq 2$ at infinity and so  a super-attracting fixpoint there. 
%, and have exactly one repelling fixpoint on their boundary (either $\infty$ or a pole of $f$).
Moreover, Lemma \ref{lemww6}  again implies that 
 all zeros of $f$ are simple.
 Assume that $W^+$ has $p$  components, all necessarily unbounded  by Lemma \ref{lemww8aa},  
 and $W^-$ has $q$ unbounded components,
 while the polynomial $f/f''$ has $r$ zeros in $\C$ arising from zeros of $f$,
all of which must be simple zeros of $f$ which are not zeros of $f''$. 

 Each component $C$ of $W^+$ has at most one pole of $f$ on its boundary, by Lemma \ref{lemww8a1},
 and so precisely one, by Lemma \ref{lemww7}  and
 the Denjoy-Wolff theorem \cite{Stei} applied to the inverse function
 $F^{-1}: H^+ \to C$, coupled with the fact that $\infty$ is a super-attracting fixpoint of~$F$ (which implies in particular that 
 $F$ is not a M\"obius transformation  and $C \neq H^+$). 
Thus  $f$ has poles at precisely $p$ points,
and each is a double zero of $f/f''$.  It now follows, in light of (\ref{wwrat3a}) and (\ref{wwrat4a}), that 
\begin{equation}
|p-q| \leq 1, \quad p+q = 1+s  \geq 2, \quad 
2p + r  = 2+s  = p+q+1, \quad r  = q - p+1, \quad 0 \leq r \leq 2. 
%2p + r - s = n+1 = p+q+1, \quad r -s = q - p+1. 
\label{pqrs1a}
\end{equation}

\begin{lem}
\label{lemratA1}
There do not exist  $x_1, x_2, x_3 \in \R$ such that $x_1 < x_2 < x_3 $ and $x_1, x_3$ are poles of $f$ while $x_2$ is a zero of $f'/f$.
\end{lem}
 \textit{Proof.} 
 Assume that such a triple  $x_1, x_2, x_3$ does exist, and without loss of generality that $f$ has no poles in
 $(x_1, x_3)$. 
No zero  of $f'/f$ can lie on the 
boundary of an unbounded component  $B$ of $W^\pm$,
 by the univalence of $F$ on $B$
 and the fact that  $F(\infty) = \infty$.  
In particular, by Lemma \ref{lemww8aa}, $x_2$ 
must lie on the boundary of a bounded component of $W^-$, and hence on the boundary 
of  a maximal finite chain $D$ of bounded components $C_j$ of $W^-$ 
joined end to end as in Definition \ref{defnchain} and its notation.  Then, in view of Lemma \ref{lemchain},
 the  $C_j$ all border the same component $A$ of $W^+$, 
and $A$ is  unbounded, with exactly one pole $x_0$ 
of $f$ on $\partial A$.  On the other hand, $\partial D$ contains no poles of $f$, and so
$x_1 < a_1 <  x_2 < b_N < x_3$. 

%080323

Suppose that $x_1 \neq x_0$.  Then $f$ has no poles in $ (x_1, a_1)$, by the choice of $x_1$ and $ x_3$, 
and  $x_1 $ lies on the boundary of some component $A' \neq A$ of $W^+$.  Since $a_1$, which is a zero of $f$, lies on $\partial A$,
 there must exist at least one zero of $F'$, and so of $f$, in $(x_1, a_1)$: let $c_1$ be the nearest such zero to $a_1$. Then there must exist a zero $d_1$ of $L = f'/f$ with $c_1 < d_1 < a_1$ and 
$L'(d_1) < 0$, which forces $d_1 \in \partial W^-$, so that  $d_1$ lies on the boundary of a bounded component of $W^-$.
Because $f$ has no zeros or poles in $(c_1, a_1)$, this 
contradicts the maximality of the finite chain $D$. 

Similar reasoning if $x_3 \neq x_0$ completes the proof of the lemma. 
 \hfill$\Box$
\vspace{.1in}

\begin{lem}
\label{lemratA2}
The integer $s$ in (\ref{wwrat3a}) satisfies $s \leq 2$, and if $s=2$ then $c < 0$. 
\end{lem}
 \textit{Proof.} 
 Suppose that $s > 2$, or $s=2$ and $c > 0$. Then there exist a large positive $R$ and 
 $\theta_j$ satisfying $0 \leq \theta_1 < \theta_2 < \theta_3 \leq \pi$, with the property  that 
 $(-1)^{j+1} ( f(Re^{i \theta_j} )- 1)$ is small, real and positive.
 %, so that 
 %$$f (Re^{i \theta_1})  \in (1, + \infty) , \quad  f (Re^{i \theta_2})  \in (0, 1) ,  \quad  f (Re^{i \theta_3}) |  \in (1, + \infty) .$$  
 
 Thus
 $Re^{i \theta_2}$ lies on a level curve $\lambda_2$  in the closed upper half-plane on which $f(z)$ is real and 
 $0 < f(z) < 1 $, and following $\lambda_2$ in the direction of decreasing $f$ leads to a real zero $y_2$ of $f$,
 possibly via one or more real zeros of $f'$.
 Similarly, for $j=1, 3$, the point  $ Re^{i \theta_j} $ lies on a level curve  $\gamma_j$ 
 in the closed upper half-plane on which $f(z)$ is real and 
 $1 < f(z) < + \infty $. Follow each $\gamma_j$ in the direction of increasing $f$: then
 $\gamma_j$ must approach a real pole $x_j$ of $f$.

 Furthermore, $\gamma_1$ and $\gamma_3$ do not meet $\lambda_2$ at all, and do not meet each other
 in the open half-plane $H^+$. 
 Hence
 it must be the case that 
 %$x_1 \neq x_3$ and
 $x_1 > y_2 > x_3$. Thus $y_2$ lies in a unbounded component $U$ of the set
 $\{ z \in \C : |f(z)| < 1 \}$, which
 cannot contain a zero of $f'$, by Lemma~\ref{lemratA1}. By the Riemann-Hurwitz formula
 \cite{Stei}, or by analytic continuation of $f^{-1}$, the function
 $f$ is univalent on $U$. But this contradicts the fact that $y_2$, $\lambda_2$ and the reflection of $\lambda_2$ across 
 $\R$ must all lie in $U$. 
% But then, since $x_1 > y_2 > x_3$, the finite boundary of $U$ must contain two unbounded
 %simple curves, on each of which $f(z) $ tends to $f(\infty) =1$ as $z $ tends to infinity in each direction. Because $\arg f(z)$ is increasing as $z$ describes $\partial U$ in the positive sense, it must be the case that $\arg f(z)$ changes by at least $4 \pi$ as $z$ describes $\partial U$, and this is a contradiction. 
 \hfill$\Box$
\vspace{.1in}

If  $s=1$ then $r=1$ and $p=1$ by (\ref{pqrs1a}), and  $f/f''$ has degree $3$, and after a real linear 
re-scaling it may be assumed that $f$ has a pole at $0$, of order $n$ say, while  the remaining zero of $f'/f''$ 
lies at $1$. Thus $f$ satisfies (\ref{hypergeometric3AB}) and is, by 
Lemma \ref{lemF4},  a constant multiple of the function 
$F_4$ in conclusion (viii) of the theorem. 

Now suppose $s=2$ and $c < 0$. Then $q=2$, $p=1$ and $r=2$ by (\ref{wwrat3a}) and 
(\ref{pqrs1a}). Hence $f/f''$ has degree $4$, with one double zero at the unique pole of $f$ and two simple zeros. 
After a linear re-scaling it may be assumed that $0$ is a pole of $f$ of order $n$, and that the simple zeros of $f/f''$ are $1$ and $K \neq 0, 1$.  This leads
 to (\ref{hypergeometric3A}) and, in view of Lemma \ref{lemF3}, to conclusion (vii),  
 which completes the proof of Proposition \ref{proprational}. 
\hfill$\Box$
\vspace{.1in}

\section{Continuation of the proof in the transcendental case} 

Assume henceforth that $f$ is transcendental and satisfies  the hypotheses of Theorem \ref{thm1}.
Since all zeros and poles of $f$ and $f''$ are real, the Tsuji characteristic of $L = f'/f$ satisfies \cite{BEL}
\begin{equation}
\mathfrak{T}(r, L ) = O( \log r ) \quad \hbox{as $r \to + \infty $.} 
\label{ww3}
\end{equation}

\begin{lem}
\label{lemww2}
The function $ f/f''$ is real entire and 
has  order of growth at most $1$. If $f/f''$ is a polynomial then $f$ satisfies
conclusion (i), (ii) or (iv) of Theorem \ref{thm1}. 
%there exist $a, b \in \R$ such that
%$f(az+b)$ is a constant multiple of $\sin z$ or $e^z$ or the function $y_1$ from (\ref{besselsoln}). 
\end{lem}
\textit{Proof.} The growth estimate follows from (\ref{ww3}) and  Lemma \ref{lemlevinost}. 
 If $f/f''$ is a polynomial then $f$ has finitely many poles and a standard application
of the Wiman-Valiron theory \cite{Hay5}
implies that $f/f''$ has degree at most $1$. If $f/f''$ is constant then evidently $f$ satisfies conclusion (i) or (ii), 
whereas if $f/f''$ has degree $1$ a real linear change of variables leads to $f$ satisfying 
equation (\ref{bessel}), in which case Lemma \ref{lembesselsoln} delivers conclusion (iv). 
 \hfill$\Box$
\vspace{.1in}

Assume henceforth that $f$ and  $f/f''$ are both transcendental. Then Lemmas \ref{lemrealzeros} and \ref{lemww2}
together imply that
\begin{equation}
\lim_{y \to +\infty, y \in \R} 
\frac{\log |f(iy)/f''(iy)|}{\log y} = + \infty . 
\label{f/f''big}
\end{equation} 

 \begin{lem}
\label{lemww2b1}
$f$ has infinitely many zeros. 
\end{lem}
\textit{Proof.} Suppose  that $f$ has finitely many zeros. Then so has $f''$, and hence 
$f'/f$ is a rational function, by the main result of \cite{La5}, and so is $f''/f$,  contrary to the assumption just made.
\hfill$\Box$
\vspace{.1in}

 \begin{lem}
\label{lemww2b}
$f$ has infinitely many poles. 
\end{lem}
\textit{Proof.} Suppose  that $f$ has finitely many poles. Then Section \ref{levostfact}  shows that $L = f'/f$ has a representation
 $L = P\psi$, where $P$ and $\psi$ are real meromorphic functions such that 
$\psi(H^+) \subseteq H^+$ and $P$ 
has finitely many poles.  Combining  (\ref{C1}) with (\ref{ww3})  delivers
$\mathfrak{T}(r, P) = O( \log r )$ as $r \to + \infty$, and so $P$ has  order at most $1$, by Lemma \ref{lemlevinost}. 

Suppose first that $P$ is transcendental with infinitely many zeros.
 Then Lemma \ref{lemrealzeros} implies that 
 $P(z)$ tends to infinity as $z \to \infty$  on $i\R^+$, faster than any power of $|z|$, 
 and so does $L(z)$, by (\ref{C1}). 
The fact that $P$ has real zeros implies that, again as $z \to \infty$ on $i\R^+$,
$$
\frac{L'(z) }{L(z)} = \frac{P'(z) }{P(z)} + \frac{\psi'(z) }{\psi (z)} = O(|z|) ,
\quad \frac{f''(z)}{f(z)} = L(z)^2 + L'(z) = L(z)^2 + O(|z|) L(z) \to \infty, 
 $$
which contradicts (\ref{f/f''big}). 

Next, if $P$ is transcendental with finitely many zeros then $zL(z)$ tends to $0$ on one of the rays 
$\arg z = \pi /4, 3 \pi /4$ and to $\infty$ on the other, by (\ref{C1}), contradicting Lemma \ref{lemnorfam}. 

Hence $P$ must be a rational function, and $f$ has finite order \cite{BEL}. Moreover,   it follows from (\ref{C1})
and (\ref{f/f''big}) that, as $z \to \infty$ on $i \R^+$, 
 $$
 L(z)^2 + L'(z) = L(z)^2 + O\left( \frac1{|z|} \right) L(z) =   \frac{f''(z)}{f(z)} = O\left( \frac1{|z|^2} \right), \quad zL(z) = O\left(1 \right) .
$$
Let $\delta $ be small and positive. Then Lemma \ref{lemnorfam} implies that, 
as $z \to \infty$ with $\delta < \arg z < \pi - \delta $, $zL(z) $ is bounded and 
$\log |f(z)| = O( \log |z| )$.  An application of the 
Phragm\'en-Lindel\"of principle now  shows that $f$ is a rational function, contrary to assumption. 
\hfill$\Box$
\vspace{.1in}

%080323
%SO I don't think we need Lemmas \ref{lemww2a}, \ref{lemww2aa}. 

%\section{Asymptotic values of $F$} 

 \begin{lem}
\label{lemww3}
The following statements hold for  asymptotic values $\beta \in \C \cup \{ \infty \}$ of $F$,
that is, values $\beta$ such that $F(z) \to \beta$ as $z \to \infty$ on a path $\Gamma_\beta$. 
\\
(i) 
There exist at most two $\beta \in \C \cup \{ \infty \}$ such that $\Gamma_\beta \cap \R$ is unbounded. \\
(ii) There exist  finitely many $\beta \in \C$ for which  $\Gamma_\beta \cap \R$ is bounded, and 
$F$ has finitely many asymptotic values.\\
(iii) All transcendental singularities of $F^{-1}$ over finite values 
are logarithmic.\\
(iv) $F$ has at most one asymptotic value $\beta \in \C \setminus \R$ 
with $\Gamma_\beta \setminus H^+$ bounded. 
\end{lem}
\textit{Proof.} To prove (i) just note that if $\Gamma_\beta \cap \R$ is unbounded
then $\beta \in \R \cup \{ \infty \}$ and it may be assumed that $\Gamma_\beta$ lies in the closed upper half plane;
hence  there is at most one $\beta$ such that 
$\Gamma_\beta \cap \R^+$ is unbounded, and at most one for which
$\Gamma_\beta \cap \R^-$ is unbounded.
Next, the first assertion of (ii) follows from
Lemma \ref{directsinglem2}, and on combination with (i) shows that 
$F$ has finitely many asymptotic values.
Since all critical points of $F$ are fixpoints of $F$, all finite singular values of $F^{-1}$ are isolated, 
so that (iii) is a consequence of the argument from
\cite[p.287]{Nev}. The fact that $F^{-1}$ has at most one direct singularity lying in $H^+$, by  Lemma \ref{directsinglem2}, then delivers (iv). 
\hfill$\Box$
\vspace{.1in}

%This leaves the question of asymptotic values $\beta  \in \R \cup \{ \infty \}$. 

\begin{lem}
\label{lemww4}
Let  $D$ be  a neighbourhood of a logarithmic singularity of $F^{-1}$ 
over $\beta \in \R$, such that $D \cap \R^+$ is unbounded. Then
there exists $a \in \R$ with $[a, + \infty ) \subseteq D$, and $f$ has finitely many zeros and poles 
on $\R^+$. 
Moreover, there cannot exist 
a neighbourhood $E \subseteq \C \setminus \R$ 
of a transcendental singularity of $F^{-1}$ 
over a finite value  $\gamma \neq \beta $.  
\end{lem} 
\textit{Proof.}
The first two assertions hold since 
$D$ is simply connected and symmetric with respect to $\R$, while all zeros and poles of $f$ are fixpoints of $F$. 

Next, assume that $E$ and $\gamma$ do exist, without loss of generality 
with $E \subseteq H^+$. 
There must exist  a path tending to infinity in $D \cap H^+$ on which $F(z) \to \beta$, and so
$F^{-1}$ has  a direct singularity over $\gamma$, lying in $H^+$, by Lemma \ref{lemww3}, 
 plus one over $\infty$,
which contradicts 
Lemma \ref{directsinglem2}. 
\hfill$\Box$
\vspace{.1in}

\begin{lem}
\label{lemww5}
The finite asymptotic values of $F$ comprise either a pair $\beta,  \overline \beta$,
where $\beta \in \C \setminus \R$, 
or one value $ \beta \in \R$.
Furthermore, all but finitely many zeros of $f$ are simple. 
\end{lem} 
\textit{Proof.} 
Suppose that 
$\beta, \gamma \in \C $ are distinct asymptotic values of $F$: then there exist   simply connected
neighbourhoods $D, E$ of 
logarithmic singularities of $F^{-1}$ 
over $\beta, \gamma $ respectively, by Lemma \ref{lemww3}.  If $D \cap \R^+$ and $E \cap R^-$ are both unbounded 
then  $\beta, \gamma $ are real and
Lemma~\ref{lemww4}, applied to $f(z)$ and $f(-z)$, implies that $f$ has finitely many poles,  contrary to assumption. 

It may therefore be assumed that either $D$ or $E$ lies in $\C \setminus \R$, and hence that 
both do, by Lemma \ref{lemww4} again. But then it must be the case that one of $D, E$ lies in 
$H^+$ and the other in the lower half-plane $H^-$,  by Lemma \ref{directsinglem2}, and moreover that
 $\gamma = \overline \beta $. 
 
 The last assertion then follows from Lemma \ref{lemww6}. 
\hfill$\Box$
\vspace{.1in}

%UP TO HERE SHOULD SUFFICE FOR THE NON-REAL AS VALUE CASE

 \begin{lem}
\label{lemww2bb}
Suppose that $f'/f$ has finitely many zeros. Then $f$ satisfies conclusion  (c) of Theorem~\ref{thm1}.  
\end{lem}
\textit{Proof.} 
This can be deduced from \cite{Laams09} 
but the following proof is included in order to keep the account self-contained. The function $f/f'$ has finitely many poles, and so has order at most $1$ by (\ref{ww3}) and Lemma \ref{lemlevinost}. On the other hand, $f/f'$ is transcendental, by Lemma~\ref{lemww2b1}. 

Since all but finitely many zeros of $f$ are simple, 
by Lemma \ref{lemww5}, the function $f$ can be written in the form $f = f_1/f_2$, in which $f_1, f_2$ are real entire functions with real zeros and no common zeros, and $f_1$ has order at most $1$. Here each $f_j$ has infinitely many zeros, by Lemmas  \ref{lemww2b1}  and \ref{lemww2b}.
Use the Levin-Ostrovskii factorisation of $f_j'/f_j$ to write
%The first step is to prove that $f_2$ has finite order, and with this aim in mind write 
$$
\frac{f'}{f} = \frac{f_1'}{f_1} -  \frac{f_2'}{f_2} =  \phi_1 \psi_1 - \phi_2 \psi_2 ,
$$
in which $\phi_j$ and $\psi_j$ are real meromorphic, while $\psi_j(H^+) \subseteq H^+$ and 
$\phi_j$ has finitely many poles.  Since $f_1$ has finite order, (\ref{C1}) leads to
$m(r, \phi_1) = O( \log r )$ as $r \to \infty$, and so $\phi_1 $ must be a rational function. 
It then follows from (\ref{C1}), (\ref{ww3}) and standard properties of the Tsuji characteristic that
$\mathfrak{T}(r, \phi_2) = O( \log r )$ as $r \to + \infty$ and so $\phi_2$ has order at most $1$ by Lemma \ref{lemlevinost}.

Let $\delta $ be small and positive and apply  Lemma \ref{lemrealzeros} to $f/f'$. On
combination with  Lemma \ref{lemnorfam} and a standard estimate for $f_1'/f_1$ \cite{Gun2}, 
this yields, as $z \to \infty$ with $\delta < \arg z < \pi - \delta$,
\begin{equation}
\label{f'/fsmallest2} 
z L(z) \to 0, \quad   \frac{f_1'(z)}{f_1(z)} = O\left( |z| \right) 
\quad \hbox{and} \quad  \frac{f_2'(z)}{f_2(z)} = O\left( |z| \right) .
%\quad  \hbox{as $z \to \infty$ with $\delta < \arg z < \pi - \delta$}.
\end{equation}
It then follows in view of (\ref{C1})
that $\log^+ | \phi_2(z) | \leq 3 \log |z| $ as $z \to \infty$ in the same sector.  Since $\delta$ may be chosen arbitrarily small, an application of the 
Phragm\'en-Lindel\"of principle now  shows that $\phi_2$ is a rational function,
so that $f_2$ has finite order \cite{BEL} and so has~$f$.

The next step is to show that $f$ and $f''$ have, with finitely many exceptions, the same zeros. Since $f''/f$ has no zeros, and all but finitely many zeros of $f$ are simple, it suffices to show that all but finitely many zeros of $f$ are zeros of $f''$. Suppose then that $x_1 , x_2 , x_3 \in \R$ are zeros of $f$ but not of $f''$, such that $x_1 < x_2 < x_3$, while
$|x_1| $ and $|x_3|$ are large and  $x_1 x_3 > 0$. Thus $x_2$ is a  simple zero of $F'$ and lies on the boundary of a component of $A$ of $W^-$. Hence it is possible to move along the real axis, away from $x_2$, while remaining on $\partial A$. 
Since $f$ cannot have a pole on $\partial A$, by Lemma \ref{lemww7},
it follows that continuing along $\R$ in the same direction until the first
encounter with a pole or zero of $f$ gives rise to a closed interval  $I \subseteq \partial A$, its endpoints being zeros of $f$. This interval $I$ must then   contain a zero of $f'/f$,  a contradiction. 

Since $f'/f$ has finitely many zeros and all but finitely many zeros of $f$ are simple,
it now follows that  the function
$$
R = \frac{f''}{ff'}
$$
has finite order and finitely many poles. As $z \to \infty$ in 
$ \delta < \arg z < \pi - \delta$,
integration of $f'/f$ using (\ref{f'/fsmallest2}), coupled with a standard estimate for $f''/f'$  from \cite{Gun2}, yields
$$
\log^+ |R(z)| \leq \log^+ \frac1{|f(z)|} + \log^+ \left| \frac{f''(z)}{f'(z)} \right| = O( \log |z| ).
$$
The Phragm\'en-Lindel\"of principle forces $R$ to be a real rational function, and 
all so but finitely many poles of $f$ are simple. Now write $R = 2/S$ and 
$$
2ff' = Sf'', \quad 
(f^2 - Sf' + S'f)' = 2ff' - S f''  -S'f' + S'f'  + S'' f = S'' f.
$$
Hence $S''f$ is the derivative of a meromorphic function and,
since $f$ has infinitely many simple poles, by Lemma \ref{lemww2b}, the rational function
$S''$ must vanish identically and 
$f^2 - Sf' + S'f$ must be a constant $c$.
This  yields  a Riccati equation 
\begin{equation}
Sf' = f^2+S'f - c = P_2(f) = (f- A_1)(f-A_2) , \quad  A_1, A_2 \in \C, \quad - c =  A_1 A_2. 
\label{riccati}
\end{equation}
If $A_1=A_2$,  then $1/S$ is the derivative of the transcendental meromorphic function $- (f-A_1)^{-1}$,
which is obviously impossible. Assume that $A_1 \neq A_2$: then 
$A_1, A_2$ are distinct Picard values, and hence asymptotic values, of $f$ and so 
$A_2 = \overline A_1$ by Lemma \ref{lemww5}. Moreover, partial fractions yields
$$
\frac{f'}{f-A_1} - \frac{f'}{f-A_2} = \frac{A_1-A_2}{ S} .
$$
Thus $S$ must be constant, since otherwise 
 $f$ is rational, a contradiction. It now follows from (\ref{riccati})  that 
$Sf' = f^2 - c$ and so $A_1+A_2 = 0$.
%$$S^2 f'' =  (2f - (A_1+\overline A_1)) (y - A_1)(y-A_2), $$
%in which $(A_1+\overline A_1)/2$ is not a Picard value of $f$. 
%Since $f''/f$ has no zeros it must be the case that $A_1+\overline A_1 =0$ 
Hence $A_1$ and $A_2 = \overline A_1$ are purely imaginary, while $- c = A_1 A_2 > 0$,
and conclusion (iii) of Theorem \ref{thm1} follows easily. 
\hfill$\Box$
\vspace{.1in}

It may be assumed henceforth that

\begin{equation}
\label{assumption}
\hbox{ $f$ and $f'/f$ each have infinitely many zeros and infinitely many poles. }
\end{equation}
 
 %090323
 
\section{Non-real  asymptotic values  } 

%This section will be occupied with the proof of the following. 

\begin{prop}
\label{propnononreal}
$F$ has no  finite non-real asymptotic values.
\end{prop}

To prove Proposition \ref{propnononreal}, assume for the remainder of this section
that $F$ has  an  asymptotic value $\beta \in \C \setminus \R$. Then by Lemma \ref{lemww5} and
the argument of 
\cite[p.287]{Nev}, it may be assumed 
that the only finite asymptotic values of $F$ are $\beta$ and $\overline \beta$ and that
there exists
an unbounded component $A$ of $W^\pm$ which contains no $\beta$-points of $F$, and which  is mapped ``infinite to one'' by $F$ onto $H^\pm \setminus \{ \beta \}$, where $H^-$ denotes the open lower half-plane. 
Moreover, the corresponding transcendental  singularity of $F^{-1}$ over $\beta$ is
logarithmic, and $A$ is simply connected. The first main step in the proof of Proposition \ref{propnononreal} will be accomplished via the following lemma.

\begin{lem}
\label{lemww10}
 There does not exist a component $B$ of $W^\pm $ such that 
 %$F \to \beta \in \C \setminus \R$ on a path tending to infinity in $B$ while
  $\pm F$ maps $B$  univalently onto $H^+$ 
 and $F(z) \to \infty$ as $z \to \infty$ on a  path in $B$.
 \end{lem}
 \textit{Proof.} Assume the contrary, let $\delta$ be small and positive, and set
 $$
 u(z) = \log^+ \frac\delta{|F(z)- \beta |} \quad (z \in A), \quad 
 u(z) = 0 \quad (z \not \in A),
 $$
 as well as $ B(t, u) = \max \{ u(z): |z| = t \}.$
 Then $u$ is subharmonic and non-constant in the plane and Lemma \ref{lemlevinost} yields, for $R \geq 1$,  
 \begin{eqnarray*}
 \int_R^{ +\infty} \frac{B(r/2, u) }{r^3} \, dr &\leq& 
  \frac3{2 \pi} \int_R^{ +\infty} \left( \int_0^{\pi} u(re^{i \theta } ) \, d \theta \right) \,  \frac{  dr }{r^3}  \\
   &\leq& 
\frac3{2 \pi}  \int_R^{ +\infty}  \left( \int_0^{\pi} \log^+ 1/|F(re^{i \theta } ) - \beta |  \, d \theta \right)
\, \frac{  dr }{r^3}  \\
& \leq&
 3 \int_R^{ +\infty} \frac{ \mathfrak{m} (r, 1/(F - \beta ) )   }{r^2} \, dr  \\ 
  &\leq& 3 \int_R^{ +\infty} \frac{ \mathfrak{T} (r, F)  + O( \log r )  }{r^2} \, dr  
  %+ O\left( \frac1{R^2} \right)   \\
     \leq  O\left( \frac{\log R}{R} \right)       .
  \end{eqnarray*}
Since $B(r/2, u)$ is non-decreasing this yields $B(R, u) = O( R \log R )$ 
as $R \to + \infty$. 

Let $\delta $ and $1/r_0$ be small and positive and denote by $\theta_A(r), \theta_B(r)$ the angular measure of the intersection with the circle $|z| = r  \geq r_0$ of $A, B$ respectively.
% (Or maybe need the set $|F-\beta| < \delta $, or even $|z^N(F-\beta)| < 1 $?\\
 %Let  $\varepsilon $ be small and positive and assume 
  Suppose first 
  that $\theta_A(r) < \pi (1 - \delta )  $ on a set $F_1$ of upper logarithmic density 
  at least $\delta  $. 
 Then, since $A \subseteq H^+$, 
 all sufficiently  large $r \in F_1$ satisfy \cite[Lemma 2.1]{BEL}
  \begin{eqnarray*}
(1+o(1)) \log r &\geq& \log B(2r, u) 
% \log m_{0\pi} (4r, 1/(F-\beta )) 
\geq 
 \int_{r_0}^r \frac{\pi dt}{t \theta_A (t)} - O(1) \\
 &\geq & 
 \int_{[r_0,r] \cap F_1} \frac{ dt}{(1- \delta ) t  } +
  \int_{[r_0,r] \setminus F_1} \frac{ dt}{t  } 
   - O(1) \\
    &\geq & 
 \int_{[r_0,r] \cap F_1} \frac{ \delta \, dt}{ (1- \delta ) t  } +
  \log r    - O(1) 
 \geq \left( \frac{\delta^2}{2 (1- \delta ) } + 1 \right)   \log r ,
 \end{eqnarray*}
 an evident contradiction.

 %It follows that 
 Hence there exists a set $E_1$ 
 of lower logarithmic density at least $1 - \delta $ on which 
 $\theta_A(r) \geq \pi (1 - \delta ) $ and so  $\theta_B (r) \leq \pi \delta  $, since $A, B$ are evidently not the same component of $W^\pm$.
 The function $w = \pm i F(z)$ maps $B$ conformally onto the right half-plane:  let $z = G(w)$ be the inverse mapping, and let 
 $\gamma_0$ be the image in $B$ under $G$
 of the real interval $[1, + \infty)$, starting from $z_0 = G(1)$. Then $\gamma_0$ tends either to infinity or to a pole of $F$, and so to infinity since $F$ is univalent on $B$. Let $r^* = |z_0| $
 and let $z = G(X) \in \gamma_0 $, with $X \geq 1$ and $r = |z|$ large. 
 Then applying Koebe's quarter theorem  to $G$ on the disc of centre $w \in [1, X] $ and radius $w$
 leads to
 \begin{eqnarray*}
 \log |F(z)| &=& \log X = \int_{[1, X]} \, \frac{|dw|}{|w|} 
 = \int_{z_0}^z \, \frac{|dz|}{ |w| |G'(w)| } \\
 &\geq& \int_{z_0}^z \, \frac{|dz|}{ 4 |z| \theta_B(|z|) } 
 \geq \int_{r^*}^{r} \, \frac{dt}{ 4t \theta_B(t)} \\
% &\geq& \int_{[r^*, r] \cap E_1} \, \frac{dt}{ 4t \theta_B(t)} \\
  &\geq& \int_{[r^*, r] \cap E_1} \, \frac{dt}{ 4 \pi \delta t } 
  \geq \left( \frac{1- 2 \delta }{4 \pi \delta } \right)  \log r  \geq 2 \log r ,
    \end{eqnarray*}
 in which
 the second integral is from $z_0$ to $z$ along $\gamma_0$.  This delivers, as $z \to \infty$ on  $\gamma_0$, 
$$|F(z) | = |z(1-1/zL(z))|  \geq |z|^{2} , \quad zL(z) \to 0.$$
On the other hand, there evidently exists a  path $\gamma_1$, tending to infinity in $A$, on which
$F(z) \to \beta $ and hence $zL(z) \to 1$. Since $L$ has only real poles, the inverse of $zL(z)$ must 
have a direct singularity over $\infty$, lying in $H^+$ and separating $\gamma_0$ from $\gamma_1$. 
But $L$ has only real zeros, and so the inverse of $zL(z)$ must 
have a direct singularity over $0$, lying in $H^+$ and separating the singularity over $\infty$ from $\gamma_1$. 
This contradicts  Lemma \ref{directsinglem2}.
   \hfill$\Box$
\vspace{.1in}

%090323

Because the component $A$ of $W^\pm$  is unbounded and simply connected and $F$ has no finite real asymptotic values,
the boundary of $A$ consists of  countably many pairwise disjoint piecewise analytic simple
curves $\gamma_j$, 
each going to infinity in both directions and
 mapped by $F$ onto
 $\R$ or  $\R \cup \{ \infty \}$, and if $F(\gamma_j) = \R \cup \{ \infty \}$ then $\gamma_j$ must meet $\R$, since $F$ has only real poles.  Suppose that one of these curves, $\gamma$ say, lies wholly in $H^+$; then $\gamma$ is mapped by $F$ onto $\R$
 and forms part of the boundary of a component $A' \neq A$ of $W^\mp $. Let $z^* \in A' $:  since $F^{-1}$ cannot have two logarithmic singularities lying in $H^+$, by Lemma \ref{directsinglem2}, analytic continuation 
 of  a local branch of the inverse of $\mp F$ shows that $z^*$ lies in a component of $W^\mp$ which is mapped univalently onto $H^+$, and which must be $A'$, so that $\gamma = \partial A'$ and $F(z)$ tends to $\infty$ along a path in $A'$, 
 contradicting Lemma \ref{lemww5}. Hence each $\gamma_j$ meets $\R$, and there is only one, because if 
 $\gamma_j$ meets $\R$ then it must separate any other  $\gamma_{j'}$ from $\R$. 
 
 Thus 
 $\partial A$ consists of a single curve, which is mapped  ``infinite to one'' onto
$\R \cup \{ \infty \}$,
 and passes in each direction through infinitely many poles of $F$, all of which are real.  
 In particular, $F^{-1} ( \{ \infty \} )$ and $\partial A \cap \R$ are neither bounded above nor bounded below,
 and neither $W^+$ nor $W^-$ has any unbounded component other than $A$.
Since $W^+$ has no bounded components, by Lemma \ref{lemww8},
while $f$ has by (\ref{assumption})
 infinitely many poles, all of which lie in $\partial W^+$ by Lemma~\ref{lemww7},
% which are then real multiple points of $F$ lying in $\partial W^+ \cap \partial W^-$,
it must be the case that $A = W^+$ and $\beta \in H^+$, and all components of $W^-$ must be bounded.

% , and   following $\gamma$  in either direction
%must   lead either to infinity or to a real simple pole of $F$, a zero of $f'/f$. Indeed, Lemma \ref{lemww5},
%elementary topological considerations and the logarithmic nature of the singularity together
 %show that all these $\gamma$ are bounded and end at poles of $F$.
%at most one of these $\gamma$ tends to infinity in both directions, while all others 
%are bounded and end at poles of $F$. 
%Thus real poles of $F$ on $\partial A $ accumulate at each of $\pm \infty$.
 
%Next, if $\lambda$ is any simple path starting at $\beta$ and tending to infinity in $H^+$

 Let  $K_0 = \{ \beta + it: \, 0 \leq t < + \infty \}$. 
 Then each pole of $F$ on $\partial A$ is the starting point of a  simple curve $\Lambda$
which
tends  
 to infinity in $A$ and is mapped injectively onto $ \{ \beta + it: \, 0 < t \leq  + \infty \}$
 by $F$. There are infinitely many of these $\Lambda$ and they are pairwise
 disjoint. Moreover, at most finitely many such $\Lambda$ meet the vertical line segment 
 $I_0 = [ {\rm Re} \, \beta, \beta ]$, because otherwise ${\rm Re} \, F$ would be constant on $I_0$ and on the curves $\Lambda$,
 contradicting  the absence of non-real critical points of $F$. 
 Choose a component $I_1$ of $\R \setminus \{ {\rm Re} \, \beta \}$ which contains infinitely many zeros of $f$:
 this is possible by 
 Lemma \ref{lemww2b1}. Because  $\partial A$  passes in each direction through infinitely many real poles of $F$, 
one of the curves $\Lambda$ can be chosen to start at a pole $y_1 \in \partial A \cap I_1$ of $F$
and not meet $I_0$. If this curve is labelled 
$K_1$ then the set $H^+ \setminus K_1$ is the union of two disjoint domains $U_1, U_2$, with   $\beta \in U_1 $ and with infinitely many zeros of $f$ lying on the boundary of  $U_2$. 
 
 \begin{lem}
 \label{lemww21}
 Choose a simple path 
 $\Gamma =
 K_2$ which starts at $\beta$,  tends to infinity and lies in $U_1$, such that $K_2$ does not meet $K_0$ except at $\beta$ itself. 
 If $y_3 \in \partial A \cap 
 \partial U_2$ is a pole of $F$ with $y_3 \neq y_1$, then there exists a path $K_3 \subseteq A
  \cup \{ y_3 \}$, 
starting at $y_3$ and tending to infinity, 
 which is mapped injectively  by $F$  onto $K_2 \cup \{ \infty \} \setminus \{ \beta \}$.  
 Moreover, $K_3$ lies in $U_2 \cup \{ y_3 \}$. 
 \end{lem}
 \textit{Proof.} 
 Here $K_2$ can be constructed using the fact that $K_1$ does not meet the  line segment 
  $I_0 $: just follow $I_0$ vertically downwards from $\beta$ and then go to infinity within $U_1$, keeping  sufficiently close to the real axis to avoid $K_1$.
  
   The existence of $K_3$ follows from analytic continuation  along $K_2$ of the branch 
 of $F^{-1}$ which maps $\infty$ to $y_3$. 
 The path $K_3$ lies in $A \cup \{ y_3 \}$ 
 and meets $U_2$, but cannot meet $K_1$ because $K_3$ and $K_1$ start at different poles of $F$ and 
 $$F( K_3 \cap K_1 ) \subseteq F(K_3) \cap F(K_1) = (K_2 \cap K_0) \cup \{ \infty \} \setminus
 \{ \beta \} =  \{ \infty \} . $$
 It follows that $K_3 \subseteq U_2 \cup \{ y_3 \}$. 
    \hfill$\Box$
\vspace{.1in}

%If $x^*$ is a zero of $f$ with $|x^*|$ large, then $x^*$ is a simple zero of $f$ and so a multiple point of $F$,
%and $x^*$  lies on $ \partial W^+ = \partial A$.
Take a zero $x_1 \in \partial U_2$ of $f$ with $|x_1|$ large, which is possible by the choice of $U_2$.
Then $x_1$ is a simple zero of $f$ and  a multiple point of $F$,
and $x^*$  lies on $ \partial W^+ = \partial A$.
Moreover,  $F(z)$  
describes $\R \cup \{ \infty \}$ monotonely and ``infinite to one'' as  the curve $\partial A$  is followed in each direction:
let $y_3, y_4$ be the first poles of $F$ which are thereby reached. 
Since $|x_1|$ is large it may  be assumed  that $|y_3|$ and $|y_4|$ are large and 
 $y_3 < x_1 < y_4$, and that $y_3, y_4 \in \partial U_2$.

 Let $\Omega = 
H^+ \setminus K_2$. Then, since $F(x_1) = x_1 \in \R$, the point 
 $x_1$ lies on the boundary of  a component $C \subseteq A $
 of $F^{-1}( \Omega )$,
 and 
 $F$ maps $C$ univalently onto
 $\Omega$, by analytic continuation of $F^{-1}$ and the fact that  $\beta \in K_2$. 
 Furthermore, the parts of $\partial A$ described in reaching $y_3, y_4$ from $x_1$ belong also to $\partial C$.
In particular,  $y_3, y_4 $ both lie in $\partial U_2 \cap \partial C \cap \partial A$.   

 Lemma \ref{lemww21} gives 
paths $K_3, K_4$ with $K_j \subseteq (A \cap U_2)  \cup \{ y_j \}$,
each starting at $y_j$ and tending to infinity, 
mapped by $F$ onto $K_2 \cup \{ \infty \} \setminus \{ \beta \}$. 
 Since $K_2$ lies in $U_1$, and no path in $C$ can cross $K_3$ or $K_4$, it follows that 
 $C$ lies in $U_2$ and thus in $\Omega$.

Now start at  $y_3$, which is a simple pole of $F$ in $\partial A \cap \R$,  
and let $z$ follow $\partial A$ in each direction 
until the first encounter with a pole of $F$
or a zero or a pole of $f$: then $z$ does not leave $\R$ as this is done, since all critical points of $F$ are zeros of $f$.
Neither of the points so reached can 
be a pole of $F$, by Lemmas \ref{lemww1} and  \ref{lemww7}, since if zeros of $L$ 
are not separated by by a zero or pole of $f$ then the 
values of $L'$ at these two zeros must differ in sign.  Thus Lemma \ref{lemww8a} implies that 
that both these points must be poles of $f$, and one of them, $y_3'$ say, lies on the part of the curve 
$\partial A$ between $y_3$ and $x_1$, which also lies in $\partial C$. 
Doing the same for $y_4$ shows that
$\partial C$ contains at least two distinct poles $y_3'$, $y_4'$ of~$f$. 
But 
 this  conclusion is incompatible with  the choice $D = \Omega$
in  Lemma \ref{lemww8a1}, giving 
a contradiction and hence completing the proof of 
Proposition \ref{propnononreal}. 

\hfill$\Box$
\vspace{.1in}

\section{Completion of the proof in the transcendental case}  
%\subsection{The case of no non-real asymptotic values}

\begin{lem}
\label{lemww31a}
%$F$ is univalent on 
All components of $W^\pm$ are mapped univalently onto $H^+$ by $\pm F$, and if $x_1$ is a  zero of $f'/f$ 
with $|x_1|$ large then $L'(x_1) < 0$ and $x_1$ 
does not lie on the boundary of a component of  $W^+$.
\end{lem}
\textit{Proof.} 
This follows from Lemmas \ref{lemww7} and \ref{lemww8}, in conjunction with  
Proposition \ref{propnononreal}.

\hfill$\Box$
\vspace{.1in}

\begin{prop}
\label{lemww31}
There do not exist sequences $x_j$ of zeros of $f'/f$ and $y_j$ of poles of $f$ both tending to $+\infty$. 
\end{prop}
\textit{Proof.}
Assume the contrary: then it is possible to choose a large positive $X_0$ and enumerate all 
the zeros $x_j$ of $f'/f$ and distinct poles $y_j$ of $f$ in $(X_0, + \infty)$ as 
$$
X_0 < x_0 < x_1 < x_2 < \ldots, \quad X_0 < y_0 < y_1 < y_2 < \ldots .
$$

Let $A_j$ be the component of $W^+$ with $y_j \in \partial A_j$. 
By Lemmas \ref{lemww8a1} and \ref{lemww31a}, it may be assumed that the $A_j$ are distinct and
their boundaries contain no zeros of $f'/f$. It then follows  that each $A_j$
contains  a path 
tending to infinity on which  $F(z) \to \infty$. Hence at most finitely many of these $A_j$ also contain a path tending to infinity on which $F(z)$ tends to a finite real asymptotic value, because otherwise $F^{-1}$ 
would have at least two direct singularities over $\infty$ lying in $H^+$,  contradicting Lemma \ref{directsinglem2}. 
Thus it may  be assumed further, for each $j$, that $\infty$ is the one and only asymptotic value approached by $F$ 
along a path tending to infinity in $A_j$. 

Similar reasoning shows that it may now also be assumed that each $x_j$ lies on the boundary of a component $B_j$ of $W^-$,
these $B_j$ being distinct and mapped univalently onto $H^+$ by $- F$. 
Hence no $B_j$  contains a path tending to infinity on which $F(z) \to \infty$, 
and again Lemma \ref{directsinglem2}
%the bound on the number of direct singularities of $F^{-1}$ in $H^+$ 
implies that at most finitely many $B_j$ contain a path on which $F(z)$
tends to a finite real asymptotic value. Thus each of these $B_j$ may be assumed to be bounded. 
%if $B_j$ not bounded can take $u_k \in B_j$ with $u_k \to\infty$ and
% $F(u_k \to w $. Then $w \in \R$ but $w \not \in \partial B_j$ so continue $F^{-1}$ along line towards $w$. 

After re-labelling if necessary, 
poles $ y_1, y_2$ of $f$ and a zero $x_m$ of $f'/f$ may be chosen 
with $y_1$ large and positive and  
$ y_1 < x_m < y_2 $.
Then $x_m$ lies on the boundary of a bounded component $B_m$ of $W^-$, and 
hence on the boundary of a maximal finite chain $D$ of bounded components of $W^-$ as in Definition \ref{defnchain}
and its notation, and Lemma \ref{lemchain} applies to $D$, with $[a_1, b_N] \subseteq \partial D \cap (y_1, y_2)$. 
Let $A$ be the component of $W^+$ given by  Lemma \ref{lemchain}: then $F$ maps $A$ univalently onto $H^+$.

Suppose first that $a_1, b_N $ are both simple zeros of $F'$. Then
 there  exist 
$c_1 < a_{1}$ and $d_N > b_{N}$ with $[c_1, a_{1}] \cup [b_{N}, d_N] \subseteq \partial A$.
Continue along $\R$ leftwards from $a_{1}$ and rightwards from $b_{N}$ until the first encounter
with a zero or pole of $f$: this is possible since $y_1 < a_1 < b_N < y_2$. 
But then  conclusion (A) of Lemma \ref{lemww8a} must hold, 
which gives at least two poles of $f$ on $\partial A$, contradicting 
Lemma \ref{lemww8a1}. 

Hence Lemma \ref{lemchain}
forces some  $x^* \in \{ a_{1}, b_{N} \}$  to be a zero of $F'$ of multiplicity  $2$, lying on the boundary 
of an unbounded component $B$ of $W^-$. Starting from $x^*$, follow the real axis, in the direction away from
$[a_1, b_N]$, until the first encounter
with a zero or pole of $f$, again possible since $y_1 < a_1 < b_N < y_2$. Then the point so reached lies on $\partial B$
and must be a zero of $f$, by Lemma \ref{lemww7}. But this gives a zero of $L= f'/f$ and hence a pole of $F$ 
lying on $\partial B$, so that $B$ is one of the $B_j$ and hence bounded, a contradiction. 
\hfill$\Box$
\vspace{.1in}

\begin{lem}
\label{lemww11a}
It may be assumed that:\\
(I) $f$ has finitely many positive poles but infinitely many negative zeros;\\
(II) $f'/f$ has infinitely many positive zeros but finitely many negative zeros;\\
(III) 
there exists a large $X_1 \in (0, + \infty) $ such that the zeros of $f$ and $f'/f$ in $(X_1, + \infty)$ are simple
and  interlaced in the sense that 
if $X_1 < a < b $ and $a, b$ are zeros of $f$ then $f'/f$ has a zero in $(a, b)$, 
while if  if $X_1 < a < b $ and $a, b$ are zeros of $f'/f$ then $f$ has a zero in $(a, b)$.
\end{lem}
\textit{Proof.} It can certainly be assumed, 
by (\ref{assumption}) and an application of Proposition \ref{lemww31} to $f(z)$
and 
$f(-z)$, that $f$ has finitely many positive poles but infinitely many negative poles, while 
$f'/f$ has finitely many negative zeros and infinitely many positive zeros. 
It then follows
% from the absence of infinitely many negative zeros of $f'/f$ 
that $f$ must also have infinitely many negative zeros, which proves (I) and (II).
Together (I) and (II) imply (III),  on  combination with
 Lemmas \ref{lemww1}, \ref{lemww7}, \ref{lemww5} and \ref{lemww31a}  and the fact that if $X_1 > 0$ is large 
 and two  zeros of $f'/f$ 
in $(X_1, + \infty)$ are
not separated by a pole or zero of $f$ or $f'/f$ then one of them has $(f'/f)' > 0$.
% and so lies on the boundary of a component of $W^+$.
 \hfill$\Box$
\vspace{.1in}

It is now possible to write 
\begin{equation}
\label{ww51}
\frac{f'}{f} = P \psi ,
\end{equation}
in which:  $\psi$ is formed as Section \ref{levostfact} using  zeros $0 < u_1 < u_2 < \ldots $
of $f$ and  zeros  $v_j \in (u_j, u_{j+1})$ of $f'/f$, and $\psi$ 
satisfies 
$\psi(H^+) \subseteq H^+$; 
the function $P$ is real meromorphic,
with finitely many zeros in $\C$, and finitely  many positive poles, but infinitely many negative poles. 
 
 \begin{lem}
 \label{lemPsmall}
 The function $P$ has order of growth at most $1$ and satisfies 
 $$
 \lim_{x \to +\infty, x \in \R} 
\frac{\log |P(x)|}{\log x} = - \infty .
$$
 \end{lem}
 \textit{Proof.}
 The first assertion follows from (\ref{ww3}), (\ref{ww51}) and Lemma \ref{lemlevinost}, applied to $1/P$.  Next,
 since $P$ is transcendental, applying 
 Lemma \ref{lemrealzeros} to $1/P$ leads to
 \begin{equation}
 \lim_{y \to +\infty, y \in \R} 
\frac{\log |P(iy)|}{\log y} = - \infty, \quad 
 \lim_{y \to +\infty, y \in \R} yL(iy) = 0.
 \label{ww66}
 \end{equation}
 Now let $\delta$ be small and positive.  Then (\ref{C1}), (\ref{ww51}) and Lemma \ref{lemnorfam}   imply that 
 $zL(z) \to 0$ and $P(z) \to 0$ as $z \to \infty$ with $  \delta <
| \arg z | < \pi -  \delta  $.
Because $P$ has finite order and finitely many poles on $\R^+$, it follows from the Phragm\'en-Lindel\"of principle that 
$P(z) \to 0$ as
 $z \to \infty$ with $ 
| \arg z | \leq  \delta  $. If  $N_1 \in \N$ is large
then (\ref{ww66}) and
the Phragm\'en-Lindel\"of principle now show that $z^{N_1}P(z)$ tends to $0$ on the sector $| \arg z | \leq  \pi/2$, 
which completes the proof. 
 \hfill$\Box$
\vspace{.1in}

Since $\psi$ maps $H^+$ into itself there exists a series representation \cite{Le}
$$
\psi(z) = az + b  +  \sum_{k=1}^\infty  A_k \left( \frac1{u_k-z} - \frac1{u_k} \right)  , 
$$
in which the $u_k$ are the poles of $\psi$, all of which are positive 
and zeros of $f$, while $a, b, A_k$ are real and $A_k > 0$, 
 $\sum_{k=1}^\infty  A_k u_k^{-2} < \infty$. On combination with Lemma \ref{lemPsmall} this implies that if 
 $k$ is large  then the residue of $f'/f = P \psi$ at $u_k$ is $-P(u_k) A_k = o(u_k^{-2} ) o( u_k^2) =
 o(1)$,
an obvious contradiction. 
This completes the proof of Theorem \ref{thm1}.

 \hfill$\Box$
\noindent
J.K. Langley, Emeritus Professor,\\
 Mathematical Sciences, University of Nottingham, \\
 NG7 2RD, UK\\
james.langley@nottingham.ac.uk
\end{document}